\documentclass[10pt]{article}

\usepackage{fancyhdr}
\newtheorem{prop}{Proposition}[section] 
\newtheorem{theo}[prop]{Th\'eor\`eme}
\newtheorem{lemm}[prop]{Lemme}
\newtheorem{coro}[prop]{Corollaire} 
\newtheorem{defi}[prop]{D\'efinition}

\usepackage[left=2.5cm,right=2.5cm,top=2cm,bottom=2cm]{geometry}

\usepackage{palatino}

\usepackage[utf8]{inputenc}
\usepackage[T1]{fontenc}     
\usepackage[francais]{babel}

\usepackage{enumerate}
\usepackage{stmaryrd}
\usepackage{amssymb} 
\usepackage{amsmath}
\usepackage{mathrsfs}
\usepackage{fancyhdr}

\usepackage{graphicx}
\usepackage{float}
\usepackage[all]{xy}
\usepackage{enumerate}

\newcommand{\bdefi}{\begin{defi}}
\newcommand{\edefi}{\end{defi}}
\newcommand{\bprop}{\begin{prop}}
\newcommand{\eprop}{\end{prop}}
\newcommand{\btheo}{\begin{theo}}
\newcommand{\etheo}{\end{theo}}
\newcommand{\blemm}{\begin{lemm}}
\newcommand{\elemm}{\end{lemm}}
\newcommand{\bconj}{\begin{conj}}
\newcommand{\econj}{\end{conj}}
\newcommand{\bcoro}{\begin{coro}}
\newcommand{\ecoro}{\end{coro}}

\newcommand{\dem}{\noindent{\bf D\'emonstration. }}
\newcommand{\cqfd}{\hfill$\Box$}

\newcommand{\rem}{\noindent{\bf Remarque. }}

\newcommand{\LLL}{{\mathcal L}}

\newcommand{\C}{{\mathcal C}}

\newcommand{\F}{{\mathcal F}}

\newcommand{\maths}[1]{{\mathbb #1}}  

\newcommand{\RR}{\maths{R}}

\newcommand{\CC}{\maths{C}}

\newcommand{\SSS}{\maths{S}}

\newcommand{\ZZ}{\maths{Z}}

\newcommand{\PP}{\maths{P}}

\newcommand{\TT}{\maths{T}}

\newcommand{\II}{\maths{I}}

\newcommand{\ttt}{{\mathfrak t}}

\newcommand{\Tr}{\operatorname{Tr}}
\renewcommand{\Re}{{\operatorname{Re}}}
\renewcommand{\Im}{{\operatorname{Im}}}
\newcommand{\im}{\operatorname{im}}

\newcommand{\Ric}{\operatorname{Ric}}

\renewcommand{\tilde}{\widetilde}

\title{Décomposition solitonique des variétés toriques}
\author{\sc{Delgove} François}
\date{\today}

\begin{document}

\maketitle

\section{Introduction}

L'article fondateur sur les solitons de Kähler-Riccci est l'article d'Hamilton \cite{Hamilton} où ils apparaissent comme points fixes du flot de Kähler-Ricci. Une manière équivalente (et que nous utiliserons par la suite) de les définir est la suivante : sur une variété compacte kählerienne de Fano $M$, le couple $(X,g)$ formé d'un champ de vecteurs holomorphe $X$ et d'une métrique kählérienne $g$ est un \textit{soliton de Kähler-Ricci} s'il existe une constante $\lambda$, dite \textit{constante d'Einstein}, telle que 
\begin{equation}\label{KRS}
\Ric(\omega_g) - \lambda \omega_g = \LLL_X \omega_g,
\end{equation}
où $\Ric(\omega_g)$ est la forme de Ricci de la forme kählérienne $\omega_g$ associée à $g$ et $\LLL_X \omega_g$ est la dérivée de Lie de $\omega_g$ dans la direction du champ de vecteurs $X$. On appellera $X$ \textit{le champ de vecteur solitonique}. Terminons ce paragraphe, en remarquant immédiatement que les solitons de Kähler-Ricci sont des généralisations des métriques de Kähler-Einstein (voir le livre \cite{besse} pour plus de détails sur les métriques de Kähler-Einstein) dont on retrouve la définition en prenant $X=0$.  

Une étape importante dans l'étude des solitons de Kähler-Ricci a été faite dans les articles \cite{TZ1,TZ2} où les auteurs introduisent un nouvel invariant de Futaki et l'utilisent pour découvrir une obstruction à l'existence de solitons de Kähler-Ricci. Plus tard, dans l'article \cite{WZ}, en utilisant à nouveau l'invariant de Futaki, les auteurs montrent, grâce à la méthode de la continuité, l'existence de solitons de Kähler-Ricci sur toute variété torique compacte de Fano.

On peut alors donner des exemples de solitons de Kähler-Ricci. L'exemple le plus simple est donné par $\CC \PP^n$, qui se trouve être une variété torique pour l'action naturelle du tore réel $\TT^n$. En effet la métrique de Fubini-Study sur $\CC \PP^n$ est une métrique de Kähler-Einstein. Une question importante est de savoir s'il existe des solitons de Kähler-Ricci qui ne soient pas des métriques de Kähler-Einstein i.e. dont le champ de vecteurs solitonique est non nul. Un théorème important est alors le théorème de Matsushima.
\btheo[\cite{Mat}]
Soit $M$ une variété de Fano admettant une métrique de Kähler-Einstein. Alors la composante neutre $\operatorname{Aut}^0(M)$ du groupe des automorphismes de $M$ est un groupe réductif complexe, et le groupe des
isométries holomorphes d'une métrique de Kähler-Einstein est un sous-groupe
compact maximal de $\operatorname{Aut}^0(M)$.
\etheo 
Or, nous savons (voir l'exemple 11.56. de \cite{besse}) que la variété torique $\CC \PP^2$ éclatée en un point a un groupe d'automorphismes non-réductif. Mais cette variété étant torique, elle admet un soliton de Kähler-Ricci qui est donc non trivial.

Rappelons maintenant que le groupe $\operatorname{Aut}^0(M)$ d'une variété compacte kählerienne $M$ est un groupe de Lie dont l'algèbre de Lie est l'ensemble des champs de vecteurs holomorphes réels que l'on notera $\eta^\RR(M)$ dans cet article. Dans le cas où la variété $M$ admet une métrique de Kähler-Einstein, le théorème de Matsushima implique donc que $\eta^\RR(M)$ est une algèbre de Lie réductive et donc la complexification de l'algèbre de Lie du sous-groupe des isométries holomorphes de $M$, or cette dernière est égale à l'ensemble des champs de vecteurs de Killing de $M$. Ce résultat a été complété par le résultat suivant (voir le lemme 28 de \cite{Pali1} ou la proposition 7.2.4 de \cite{dercalabi}). Avant de l'énoncer, nous avons besoin de normaliser l'équation des solitons de Kähler-Ricci \eqref{KRS} en prenant $\lambda=1$ et d'introduire la notion de laplacien pondéré. On définit le \textit{laplacien pondéré} par $X \in \eta^\RR(M)$ grâce à la formule suivante pour $u \in \C^\infty(M,\RR)$
$$
\Delta^X_{g,J} u:= \Delta_g u - \cfrac{1}{2} \, X^{1,0}(u). 
$$
De plus, puisque la structure complexe $J$ de la variété $M$ est antisymétrique pour $g$, on peut étendre la définition précédente par $\CC$-linéarité aux fonctions $u \in \C^\infty(M,\CC)$. Nous avons alors une décomposition de l'espace des champs de vecteurs holomorphes réels que nous appellerons par la suite décomposition solitonique. On pourra consulter le théorème \ref{decomposolito} pour plus de détails. 

L'objectif principal de cet article est d'étudier cette décomposition dans le cas torique. Rappelons ici que $(M,\omega, J, \TT^n, \mu)$ est une variété torique kählérienne (connexe) de dimension $2n$ si
$(M,\omega)$ est une variété kählérienne (connexe) compacte de dimension $2n$ munie d'une action hamiltonienne et fidèle du tore $\TT^n$ de dimension $n$ dont on note $\ttt$ l'algèbre de Lie et d'application moment $\mu : M \rightarrow \ttt^*$. Dans le cas torique, voir par exemple \cite{WZ,Donaldson}, rechercher un soliton de Kähler-Ricci $(X,g)$ revient à chercher les couples $(g,a)$ composés d'une métrique kählérienne $g$ sur $M$ et d'un vecteur $a \in \ttt$ vérifiant
\begin{equation}
Scal_g - \overline{Scal} = -2 \, \Delta^g \langle \mu, a \rangle,
\end{equation}
où $Scal_g$ est la courbure scalaire de $g$ et $\overline{Scal}$ la courbure scalaire totale. On notera alors $\Delta^{g,a}$ le laplacien pondéré et $\Delta^a_{g,J}$ le laplacien complexe pondéré (voir les équations \eqref{lapcomgen} et \eqref{lapcomactang}). Avec cette notation, nous montrerons le lemme suivant:
\blemm
Soit une variété torique kählérienne compacte de Fano $(M,\omega_0,J,\TT^n, \mu)$. Alors le couple $(g,a)$ est un soliton de Kähler-Ricci si et seulement si modulo une constante  additive on a
$$
\forall b \in \ttt,~~ 2 \langle  \mu, b \rangle = \Delta^{g,a} \langle \mu, b \rangle,
$$
\elemm

Rappelons que l'image de l'application moment d'une variété torique est un polytope de Delzant et que de plus, quitte à ajouter une constante au polytope moment, on peut supposer que ce polytope est égal à l'opposé du polytope induit par son éventail si on considère la variété comme une variété algébrique torique (voir section \ref{constructionagl}). On peut alors considérer les racines de Demazure de ce polytope, ces dernières vont servir à paramétrer les espaces propres du laplacien pondéré complexe.

\btheo
Soit $(M,\omega_0,J,\TT^n,\mu)$ une variété torique kählérienne compacte de Fano telle que $P:=\im(\mu)$ soit le polytope algébrique associé, dont on note $R(P)$ l'ensemble de ses racines de Demazure. Supposons que $(g,a)$ soit un soliton de Kähler-Ricci de constante de Kähler-Einstein égale à $1$ sur la variété complexe $(M,J)$. 
\begin{itemize}
\item[$\bullet$] Nous avons la décomposition suivante :
$$
\overline{\ker \left( \Delta^a_{g,J} - 2 \II \right)} = \operatorname{Aff}_0^\CC \oplus \bigoplus_{\alpha \in R(P)} \CC \, \tilde{v_\alpha},
$$
où $\operatorname{Aff}_0^\CC$ est l'ensemble des fonctions lisses sur $M$ s'écrivant sous la forme $\langle x , b_1 \rangle + \sqrt{-1} \langle x , b_2 \rangle$ où $(b_1,b_2) \in \ttt \times \ttt$ et $\tilde{v_\alpha}$ est la fonction appartenant à $\C^\infty(M,\CC)$ telle qu'en coordonnées action-angles (x,t)
$$
\tilde{v_\alpha} \vert_{M^0} = \left( \langle x, b_{\rho_\alpha} \rangle + 1 \right) e^{- \langle \alpha , \nabla \phi - \sqrt{-1} \, t \rangle },
$$
où $\phi$ est le potentiel symplectique associé à $g$.

\item[$\bullet$] Si on note $V_{\gamma_i}$ les espaces propres de la décomposition solitonique et $\chi$ l'isomorphisme entre $ \overline{\ker \left( \Delta^a_{g,J} - 2 \II \right)}$ et $\eta^{\RR}(M)$ du théorème \ref{decomposolito}, alors nous avons
$$
\chi^{-1}(V_0)= \operatorname{Aff}_0^\CC \oplus \bigoplus_{\alpha \in R(P), \langle \alpha, a \rangle=0} \CC \, \tilde{v_\alpha}
$$
et pour tout $i \geq 1$,
$$
\chi^{-1}(V_{\gamma_i})= \bigoplus_{\alpha \in R(P), \, 2 \langle \alpha, a \rangle=\gamma_i} \CC \, \tilde{v_\alpha}.
$$
\item[$\bullet$] Les fonctions $\tilde{v_\alpha}$ sont des fonctions propres de l'opérateur $\Delta^a_{g,-J} - 2 \II $ pour la valeur propre $4 \langle a , a \rangle$.
\end{itemize}
\etheo

Nous terminons cet article par l'étude de deux exemples : l'espace projectif $\CC \PP^2$ et son éclaté en un point.

\textit{Remerciements.} Cet article est extrait de la thèse de l'auteur réalisé sous la direction de Nefton Pali. L'auteur tient aussi à remercier Frédéric Paulin pour sa relecture attentive et ses nombreux commentaires pertinents.
\section{Rappels de géométrie torique}\label{torsymp}

\subsection{Géométrie torique kählérienne}
Toutes nos variétés sont supposées connexes sauf mention explicite du contraire. Rappelons qu'un quadruplet $(M,\omega, \TT^n, \mu)$ est une variété torique symplectique de dimension $2n$ si
$(M,\omega)$ est une variété symplectique compacte de dimension $m:=2n$ munie d'une action hamiltonienne et fidèle du tore $\TT^n$ de dimension $n$ et d'application moment $\mu$. On note $\ttt$ l'algèbre de Lie de $\TT^n$ et $\Lambda$ le réseau déterminé par le tore $\TT^n$ i.e. $\TT^n=\ttt/\Lambda$. Dire que $\mu : M \rightarrow \ttt^*$ est une application moment signifie qu'elle est $\TT^n$-invariante et qu'elle vérifie 
\begin{equation}\label{fanfan1}
\forall b \in \ttt,~~-d \langle \mu , b \rangle = i_{X_b} \omega,
\end{equation}
où 
$X_b(z):= \frac{d}{dt} \vert_{t=0} \exp (t b) \cdot z$ est le champ de vecteurs fondamental de l'action du tore $\TT^n$ associé au vecteur $b \in \ttt$. Précisons aussi qu'elle est unique à une constante additive près.

Un \textit{polytope de Delzant} est un polytope convexe dans $\ttt^*$ de $\ttt$ défini par $d>n$ inégalités de la forme
$$
\langle \nu_r , x \rangle + \lambda_r \geq 0,
$$
où $\nu_r$ appartient à $\Lambda$ et $\lambda_r \in \RR$ et tel que tout sommet soit l'intersection de $n$ faces de codimension 1 dont les vecteurs normaux forment une base de $\Lambda$.
\bprop[voir par exemple \cite{audin2012torus}]
Nous avons les propriétés suivantes:

\begin{enumerate}\label{proptor}
\item[$\bullet$] L'image de $\mu$ dans $\ttt^*$ est un polytope convexe $P$ qui est l'enveloppe convexe des images par $\mu$ des points fixes de l'action de $\TT^n$. On montre que ce polytope est un polytope de Delzant et que réciproquement tout polytope de Delzant est l'image par une application moment d'une variété torique.

\item[$\bullet$] Pour chaque face $F$ de $P$, si on note $F^0$ l'intérieur de $F$ alors le stabilisateur de tout point $z \in \mu^{-1}(F^0)$ pour l'action de $\TT^n$ est le tore de dimension égale à la codimension de $F$ et dont l'algèbre  de Lie est l’annulateur dans $\ttt$ du sous-espace vectoriel de $\ttt^*$ correspondant à $F$.

\item[$\bullet$] En particulier, si on note $P^0$ l'intérieur de $P$, alors l'action de $\TT$ sur $M^0:=\mu^{-1}(P^0)$ est libre et donc $\mu : M^0 \rightarrow P^0$ est une $\TT$-fibration principale. De plus, $M^0$ est un ouvert dense de $M$.
\end{enumerate}
\eprop

Ainsi si $P$ est un polytope de Delzant, alors il s'écrit sous la forme 
\begin{equation}\label{Polytopgen}
P = \bigcap_{i=1}^d \lbrace x \in \ttt^* ~:~ L_i(x):=\langle x, \nu_i  \rangle + \lambda_i \geq 0 \rbrace,
\end{equation}
où $d>n$, $\nu_i \in \Lambda$ et $\lambda_i \in \RR$. Une face est alors déterminée par un ensemble $I_F \subset \lbrace 1, \cdots, r \rbrace$ tel que
$$
F = P \cap \bigcap_{i \in I_F} \lbrace x \in \ttt^* ~:~ L_i(x)= 0 \rbrace.
$$
On notera $\F(P)$ l'ensemble des faces du polytope $P$. Le deuxième point de la proposition \ref{proptor} signifie que pour toute face $F \in \F(P)$ d'intérieur $F^0$, si on note $\TT_F := \ttt_F /\Lambda_F$ où $\ttt_F $ est le sous-espace vectoriel de $\ttt$ engendré par les $\lbrace \nu_i ~:~ i \in I_F \rbrace$ et $\Lambda_F$ est le sous-réseau de $\Lambda$ engendré par les $\lbrace \nu_i ~:~ i \in I_F \rbrace$, alors
$$
\mu^{-1}(F^0)= F^0 \times \TT^n/\TT_F.
$$
Le troisième point de la proposition \ref{proptor} signifie que $M^0$ est symplectomorphe à $P^0 \times \TT^n \subset \RR^n \times \TT^n$ que l'on a muni de la forme symplectique induite par la forme symplectique standard de $\RR^{2n}$. Ainsi en identifiant $\ttt \simeq \RR^n$ et $\Lambda \simeq \ZZ^n$, nous obtenons que
$$
M^0 \simeq P^0 \times \TT^n = \left\lbrace (x,t) \in P^0 \times \RR^n/\ZZ^n \right\rbrace.
$$
Ce système de coordonnées sur $M^0$ est appelé le système de \textit{coordonnées actions-angles}. L'action de $\TT^n$ est alors donnée par
$$
\forall \theta \in \TT^n,~~\theta \cdot (x,t) = (x, t+ \theta),
$$
et la forme symplectique $\omega$ s'écrit simplement (en restriction à $M^0$) sous la forme
$$
\omega = \sum_{j=1}^n dx_j \wedge dt_j.
$$
On peut l'écrire matriciellement sous la forme :
\begin{equation}\label{omegamatrix}
\omega
=
\left(
\begin{array}{cc}
0 & I \\ 
-I & 0
\end{array} 
\right).
\end{equation}
Remarquons que, dans le système de coordonnées actions-angles, nous avons :
\begin{equation}
\forall b \in \ttt,~~ X_b = \sum_{i=1}^n b_i \cfrac{\partial}{\partial t_i}.
\end{equation}

On rajoute une structure complexe $J$ qui soit $\TT^n$-invariante et telle que $(M,\omega,J)$ soit une variété kählérienne i.e. $g:=\omega( \cdot, J \cdot)$ est une métrique riemannienne. On dira alors que \textit{$(M,\omega,J,\TT^n,\mu)$ est une variété torique kählérienne}. 

La première conséquence est que, puisque $J$ est intégrable, l'action du tore $\TT^n$ s'étend en une action holomorphe du tore complexe $\TT_\CC \simeq (\CC^\times)^n$. En particulier, nous obtenons un système de coordonnées holomorphes $(u,v)$ sur $M^0$:
\begin{equation}\label{fanfan2}
M^0 \simeq (\CC^\times)^n \simeq \RR^n \times \sqrt{-1} \, \TT^n = \left\lbrace u + \sqrt{-1} \, v ~:~  u \in \RR^n,~~ v \in \RR^n/\ZZ^n \right\rbrace.
\end{equation}
Dans ce système, nous avons
$$
\forall \theta \in \TT^n,~~\theta \cdot \left( u + \sqrt{-1} v \right) = u + \sqrt{-1} \, \left( v + \theta \right).
$$
La structure complexe $J$ est alors simplement donnée par la multiplication par $\sqrt{-1}$ i.e.
\begin{equation}
J
=
\left(
\begin{array}{cc}
0 & -I \\ 
I & 0
\end{array} 
\right).
\end{equation}
De plus, grâce à l'équation \eqref{fanfan2}, on sait que $M^0 \simeq \TT_\CC$, cette dernière étant une variété de Stein, il existe un potentiel  $f \in \C^\infty(M^0,\RR)$ tel que 
\begin{equation}\label{potf}
\omega = 2 \sqrt{-1} \partial \overline{\partial} f.
\end{equation}
Puisque $\omega$ est $\TT^n$-invariante, le potentiel $f$ ne va dépendre que de la variable $u$. Ainsi nous obtenons
\begin{equation}
\omega
=
\left(
\begin{array}{cc}
0 & F \\ 
-F & 0
\end{array} 
\right)\text{ et }
g 
=
\left(
\begin{array}{cc}
F & 0 \\ 
0& F
\end{array} 
\right).
\end{equation}
où $F$ est la matrice hessienne de $f$ pour la variable $u$.
De plus, nous allons travailler sur des variétés de Fano, cette propriété est équivalente (voir par exemple \cite{Donaldson,WZ}) sur le polytope moment $P$ à l'existence d'un \textit{centre privilégié} : c'est-à-dire avec les notations de l'équation \eqref{Polytopgen} qu'il existe $x \in \ttt^*$ tel que $L_1(x)= \cdots = L_d(x)$. Un tel $x$ est alors unique. On définit alors \textit{l'ensemble des potentiels symplectiques} $S(P)$ comme l'ensemble des fonctions $\phi \in \C^{0}(P)$ telles que sa restriction à $P^{0}$ ou à l'intérieur de toute face non vide de $P$ soit lisse et strictement convexe et telle que si
$\phi_0= \frac{1}{2} \, \sum_{k=1}^d L_k \log L_k$ alors $\phi - \phi_0$ soit la restriction à $P$ d'une fonction lisse définie sur un ouvert contenant $P$. Nous avons le résultat suivant (on pourra consulter \cite{Abreu}).

\btheo\label{potsymp}
L'ensemble des métriques $\TT^n$-invariantes  sur $(M,\omega,\TT^n,\mu)$ est en bijection avec le quotient de $S(P)$ par l'ensemble $\operatorname{Aff}(P,\RR)$ des sur $P$. De plus la bijection est donnée par l'application $\phi \mapsto g_\phi$ où
\begin{equation}\label{gtorique}
\forall \phi \in S(P),~~ g_\phi \vert_{M^0} = \sum_{i,j} G_{ij} d x_i \otimes d x_j + H_{ij} dt_i \otimes dt_j,
\end{equation} 
dans le système de coordonnées actions-angles où $G=(G_{ij})_{ 1 \leq i,j \leq j}$ est la matrice hessienne de $\phi$ et $H=(H_{ij})_{ 1 \leq i,j \leq j}$ est la matrice inverse de  la matrice $G$. 
\etheo
On dira alors que $\phi$ est le \textit{potentiel symplectique} de $g_\phi$, uniquement déterminé modulo $\operatorname{Aff}(P,\RR)$.
Commençons par remarquer que le potentiel $\phi_0$ dont l'expression est donnée par
\begin{equation}\label{Guillemin}
\phi_0= \frac{1}{2} \, \sum_{k=1}^d L_k \log L_k,
\end{equation}
définit un potentiel symplectique que l'on appelle le \textit{potentiel de Guillemin}. De plus, on peut écrire matricellement les différentes expressions précédentes dans le système de coordonnées actions-angles :
\begin{equation}\label{gmatrix}
\omega
=
\left(
\begin{array}{cc}
0 & I \\ 
-I & 0
\end{array} 
\right),~
g
=
\left(
\begin{array}{cc}
G & 0 \\ 
0 & H
\end{array} 
\right)
\text{ et }
J
=
\left(
\begin{array}{cc}
0 & -H \\ 
G & 0
\end{array} 
\right).
\end{equation}
\\
Ces expressions nous permettent de calculer les gradients riemannien et symplectique d'une fonction $\psi \in \C^\infty(M,\RR)$. Nous obtenons donc
\begin{equation}\label{nablag}
\nabla_g \psi = \sum_{i,j=1}^n H_{ij} \dfrac{\partial \psi}{ \partial x_j}  \dfrac{\partial}{ \partial x_i} + \sum_{i,j=1}^n G_{ij} \dfrac{\partial \psi}{ \partial t_j} \dfrac{\partial}{ \partial t_i} 
\end{equation}
et
\begin{equation}\label{nablaw}
\nabla_\omega \psi =\sum_{i=1}^n \dfrac{\partial \psi}{ \partial x_i} \dfrac{\partial}{\partial t_i} - \sum_{i=1}^n \dfrac{\partial \psi}{ \partial t_i}  \dfrac{\partial}{\partial x_i}.
\end{equation}
On peut aussi calculer la forme de Ricci $\Ric(\omega)$ dans ce système de coordonnées, nous obtenons donc
\begin{equation}\label{Ricactang}
\Ric(\omega)= -\cfrac{1}{2} \sum_{i,j,k,l=1}^n \, H_{li,ik} \, d x_k \wedge dt_l \text{ avec }
H_{li,ik} := \cfrac{ \partial H_{li}}{\partial x_i \partial x_k}.
\end{equation}

Terminons par \textit{la formule du laplacien} pour la métrique $g$ :
\begin{equation}\label{Laplacien}
\Delta^{g}= - \sum_{i,j=1}^n \cfrac{\partial}{\partial x_i} \left( H_{ij} \cfrac{\partial}{\partial x_j} \right) + \cfrac{\partial}{\partial t_i} \left( G_{ij} \cfrac{\partial}{\partial t_j} \right),
\end{equation}
et par \textit{la formule d'Abreu} pour la courbure scalaire :
\begin{equation}\label{Abreu}
Scal_g = \mu^{*}S(H), \text{où } S(H)= - \sum_{i,j=1}^n \cfrac{\partial^2 H_{ij}}{\partial \mu_i \partial \mu_j}.
\end{equation}

Terminons cette section en faisant le lien entre les deux systèmes de coordonnées $(x,t)$ et $(u,v)$. En effet, on peut remarquer que l'application moment $\mu$ est donnée en coordonnées complexes par :
\begin{equation}\label{Paulin1}
\mu(u,v)=\dfrac{\partial f}{\partial u}.
\end{equation}
Ainsi par restriction, on peut voir que $\mu(u,0)$ est un difféomorphisme entre $\RR^n$ et $P^0$. Le changement de coordonnées est donc donné par le difféomorphisme suivant :
\begin{equation}\label{Pierre3}
x= \dfrac{\partial f}{\partial u},~~ t=v.
\end{equation}
On pourra consulter \cite{Abreu,calderbank} pour plus de détails.
\subsection{Construction algébrique}\label{constructionagl}

Cette section qui s'inspire de l'article \cite{mabuchi1987} rappelle la construction "algébrique" d'une variété torique kählérienne de Fano. Cette construction nous apporte des précisions sur l'application moment et sur la décomposition du groupe des automorphismes d'une variété torique (et de l'algèbre de Lie des champs de vecteurs holomorphes).

\subsubsection{Éventails et construction algébrique}

Soit $A$ un $\ZZ$-module libre de rang $n$. On note $B:=\operatorname{Hom}_\ZZ(A,\ZZ)$ le groupe des morphismes (de modules) de $A$ dans $\ZZ$. On  a alors un accouplement bilinéaire $\langle \cdot , \cdot \rangle : B \times A  \rightarrow \ZZ$.

On peut étendre $B$ et $A$ en deux espaces vectoriels réels de dimension $n$ en prenant leur produit tensoriel avec $\RR$: on pose donc
$ B_\RR := B \otimes_\ZZ \RR$ et $ A_\RR := A\otimes_\ZZ \RR$. L'accouplement $\ZZ$-bilinéaire s'étend alors en une application 
$\RR$-bilinéaire entre $B_\RR$ et $A_\RR$. Si $\sigma$ est sous-ensemble de $A_\RR$ alors on dit que $\sigma$ est un \textit{cône} s'il existe $(a_1,\cdots,a_s) \in A^s$ tel que 
\begin{itemize}
\item[$\bullet$] $\sigma = \RR^+ \, a_1 + \cdots + \RR^+a_s$,
\item[$\bullet$] $\sigma \cap (-\sigma) = \lbrace 0 \rbrace$.
\end{itemize}
On définit alors \textit{le cône dual} de $\sigma$ par $\sigma^\vee:=\lbrace x \in B_\RR ~:~ \langle x , y \rangle \geq 0 ~~\forall y \in \sigma \rbrace$. De plus, on dit qu'un sous-ensemble $\tau$ de $\sigma$ est \textit{une face} de $\sigma$, et on note $\tau \leq \sigma$, s'il existe $b_0 \in \sigma^\vee$ tel
$
\tau = \sigma \cap  b_0^\perp = \lbrace y \in \sigma ~:~ \langle b_0 , y \rangle = 0 \rbrace.
$

\bdefi
Un éventail de $A$ est un ensemble $\Delta$ de cônes de $A_\RR$ telle que
\begin{itemize}
\item[$\bullet$] $\forall \sigma \in \Delta, \tau \leq \sigma \Rightarrow \tau \in \Delta$
\item[$\bullet$] $\forall (\sigma,\sigma') \in \Delta^2,~~ \sigma \cap \sigma' \leq \sigma \text{ et } \sigma \cap \sigma' \leq \sigma'$.
\end{itemize}
\edefi

Pour tout éventail $\Delta$, on définit \textit{le support} $\vert \Delta \vert$  de $ \Delta$ dans $A_\RR$ par
$
\vert \Delta \vert := \cup_{\sigma \in \Delta} \sigma. 
$
De plus pour tout $i \in \lbrace 0, \cdots, n \rbrace$, on pose 
$$
\Delta(i):= \lbrace \sigma \in \Delta ~:~ \dim \sigma = i \rbrace,
$$
où $\dim \sigma$  désigne la dimension de l'espace vectoriel réel engendré par $\sigma$ dans $A_\RR$. On dira que $\Delta$ est \textit{un éventail non-singulier}  si pour tout $\sigma \in \Delta(n)$, il existe une $\ZZ$-base $(a_1 , \cdots, a_n )$ de $A$ et un entier $s \leq n$ tels que $\sigma = \RR^+ \, a_1 + \cdots + \RR^+ \, a_s$. Cela implique qu'il existe un ensemble fini $(a_1, \cdots, a_s) $ d’éléments irréductibles dans $A$ engendrant $\sigma$ en tant que cône positif. On dira que les $(a_1, \cdots, a_s)$ sont \textit{les générateurs fondamentaux} de $\sigma$.

Avec ces notations, on a alors le théorème suivant, que l'on tire du théorème 1.4 de \cite{mabuchi1987}, mais dont on peut trouver des preuves dans \cite{Demazure1970}, \cite{oda2012convex}.

\btheo
Soit $\Delta$ un éventail non singulier de $A$. On peut lui associer une unique variété algébrique complexe lisse  $\TT_\Delta$ qui vérifie les propriétés suivantes.
\begin{itemize}
\item[$\bullet$] $\TT_\Delta$ est une compactification $\TT_\CC$-équivariante de $\TT_\CC$ irréductible de dimension $n$. 
\item[$\bullet$]  Pour tout $i \in \lbrace 0, \cdots, n \rbrace$ et pour tout $\sigma \in \Delta(i)$, il existe une unique orbite $O^\sigma$ sous l'action de $\TT_\CC$ telle que
$$
\TT_\Delta = \bigsqcup_{\sigma \in \Delta} O^\sigma.
$$
De plus, l'adhérence $D(\sigma)$ de $O^\sigma$ est une sous-variété $\TT_\CC$-stable de $\TT_\Delta$ topologiquement irréductible lisse  et de dimension $n-i$ qui admet la décomposition suivante :
$$
D(\sigma) = \bigsqcup_{\tau \geq \sigma} O^\tau
$$
\item[$\bullet$] Pour tout $\sigma \in \Delta(n)$, $U_\sigma := \cup_{\tau \leq \sigma} O^\tau$ est un voisinage affine ouvert et $\TT_\CC$-stable de $O^\sigma$ dans $\TT_\Delta$ tel que 
$$
\TT_\CC \subset U_\sigma \simeq \CC^n,
$$
et
$$
\TT_\Delta = \bigcup_{\sigma \in \Delta(n)} U_\sigma.
$$
\end{itemize}
\etheo

Nous avons alors la réciproque au  théorème précédent (voir le théorème 4.1 de \cite{MO}):
\btheo
Toute variété $X$ algébrique complexe lisse irréductible  de dimension $n$, sur laquelle  $\TT_\CC$ agit de manière fidèle et régulière, est $\TT_\CC$-équivariament isomorphe à une variété de la forme $\TT_\Delta$ pour un éventail $\Delta$ non-singulier d'un $\ZZ$-module libre $A$ de rang $n$.
\etheo

Terminons cette section par une caractérisation des variétés de Fano.  

\btheo\label{Palg}
Soit  une variété algébrique complexe irréductible compacte torique $M$ lisse dont l'éventail est noté $\Delta$ i.e. $M = \TT_\Delta$ où $\Delta$ est un éventail d'un $\ZZ$-module libre $A$ de rang $n$. On a l'équivalence entre les propositions suivantes :
\begin{enumerate}
\item $M$ est une variété de Fano,
\item $P_{alg}:=  \lbrace a \in \operatorname{Hom}(A,\ZZ)_\RR ~:~ \langle a , b_\rho \rangle \leq 1 ~~\forall \rho \in \Delta(1) \rbrace$ est un polytope convexe compact dont les sommets sont les $\lbrace a_\tau ~:~ \tau \in \Delta(n) \rbrace$, où $a_\tau$ est l'unique élément de $\operatorname{Hom}(A,\ZZ)_\RR$ tel que $ \langle a_\tau,b \rangle=1$ pour tout générateur $b$ de $\tau$.
\end{enumerate}
\etheo

\subsubsection{Lien entre le polytope algébrique et le polytope symplectique}\label{secPolAlg}

Dans cette section, nous faisons, à l'aide de l'application moment, le lien avec la construction symplectique et avec le polytope de Delzant associé à cette construction.\\

Commençons par rappeler que
$
\TT_\CC
$
est un ouvert dense de $\TT_\Delta$. De plus, on sait qu'il existe un sous-groupe compact maximal $\TT_\RR$ dans $\TT_\CC$ tel que
$$
\lbrace (t_1, \cdots, t_n) ~:~ t_i \in \CC^\times, \vert t_i \vert =1 \rbrace = \SSS_1^n \simeq \TT_\RR \subset \TT_\CC \simeq (\CC^\times)^n = \lbrace (t_1, \cdots, t_n ) ~:~ t_i \in \CC^\times \rbrace
$$

Si on note $(t_1, \cdots, t_n) \in (\CC^\times)^n$ le système de coordonnées holomorphes usuelles de $\TT_\CC$, alors on peut définir, pour tout $i \in \lbrace 1, \cdots, n \rbrace$, des fonctions $x_i \in \C^\infty(\TT_\CC,\RR)$ par la formule suivante :
$$
t_i \, \overline{t_i} = \exp(-x_i).
$$
Ainsi, toute fonction $u \in \C^\infty(\RR^n,\CC)$ en les variables $(x_1,\cdots, x_n)$ peut être vue comme une fonction appartenant à $\C^\infty(\TT_\CC^n,\CC)$ en les variables complexes $(t_1, \cdots,t_n)$ qui sera $\TT_\RR$-invariante.

Considérons maintenant la fonction $u^0$ définie sur $\RR^n$ (identifié à $A_\RR$ par le choix d'une $\ZZ$-base de $A$) par
$$
u^0 : x \in \RR^n  \mapsto \log [ \displaystyle \sum_{\tau \in \Delta(n)} \exp ( \langle a_{\tau},x \rangle )] \in \RR 
$$
Un calcul direct nous permet de voir que cette fonction est strictement convexe  sur $\RR^n$ et que $ \mu_{u^0}:=\nabla u^{0}: x \in \RR^n \mapsto ( \partial_1 \, u^0 (x), \cdots, \partial_n \, u^0 (x)) \in \RR^n$ (où $\RR^n$ est identifié $ \operatorname{Hom}(A,\ZZ)_\RR$ par le choix de la $\ZZ$-base duale de la $\ZZ$-base choisie pour $A$) est un difféomorphisme entre  $\RR^n$ et l'intérieur $P_{alg}^0$ du polytope $P_{alg}$ défini au théorème \ref{Palg}. La $2$-forme $\omega_0$ sur $\TT_\CC$ définie (via l'identification précédente) par
$$
\omega_{0}:= \sqrt{-1} \, \partial \overline{\partial} \, u^0,
$$
définit une métrique kählérienne $\TT_\RR$-invariante sur $\TT_\CC$. Nous avons alors le résultat suivant :

\btheo
La métrique $\omega_{0}$ s'étend sur $M$ en une métrique kählérienne $\TT_\RR$-invariante vérifiant en plus $\omega_{0} \in 2 \pi \, c_1(M)$. De plus, on a $\im(\mu_{u^0})= -P_ {alg}$.
\etheo
\dem
On pourra consulter la section 3 de \cite{BS} pour la première assertion et le théorème 4.2 de \cite{mabuchi1987} pour la seconde assertion.
\cqfd
\medskip

De plus, si on se donne une autre métrique kählérienne  $\TT_\RR$-invariante $g$ dont la $(1,1)$-forme $\omega_g$ appartient à $ 2 \pi \, c_1(M)$ alors d'après le $\partial \overline{\partial}$-lemme, il existe $\varphi \in \C^{\infty}(M,\RR)$ telle que $\omega = \omega_{0} + \sqrt{-1} \partial \overline{\partial} \varphi$. En remarquant que $\varphi$ est donc $\TT^n$-invariant, nous avons alors, sur $\TT_\CC$, l'égalité
$$
\omega = \sqrt{-1} \, \partial \overline{\partial} u , \text{ où } u= u^0 + \varphi.
$$
On retrouve ainsi le potentiel kählérien modulo un facteur multiplicatif $2$ (voir l'équation \eqref{potf}).
De plus, puisque la fonction $\varphi$ est globalement définie sur $M$, elle est bornée et donc on a $\im(\mu_{u})= \im(\mu_{u^0}) = P_{alg}$. On pourra consulter \cite{BS} ou \cite{WZ} pour plus de détails.

On a le corollaire suivant qui fait le lien avec la construction symplectique.

\bcoro
Soit $(M,\omega,J,\TT^n,\mu)$ une variété kählérienne compacte torique de Fano, d'application moment $\mu$  telle que $\omega \in 2 \pi \, c_1(M)$.
On peut associer à cette variété deux polytopes $P_{symp}= \mu(M)$ et le polytope $P_{alg}$ provenant d'un éventail $\Delta$ construit à la section précédente. Nous avons alors
$$
P=P_{symp}= - P_{alg}.
$$
En particulier, si on écrit que
$$
P_{symp} = \bigcap_{r=1}^d \lbrace x \in \RR^n ~:~ \langle x , \nu_r \rangle \geq - \lambda_r \rbrace,
$$
alors $d$ est égal au cardinal de $\Delta(1)$ et pour tout $r=1,\cdots,d$, il existe un unique $\rho \in \Delta(1)$ tel que
$\nu_r = b_\rho \text{ et } \lambda_r=1$ i.e.
\begin{equation}\label{Polytoperef}
P_{symp} = \bigcap_{\rho \in \Delta(1) } \lbrace x \in \RR^n ~:~ \langle x , b_\rho  \rangle \geq -1 \rbrace.
\end{equation}
\ecoro

\dem
Avec les notations précédentes, puisque $t_i \overline{t_i} = \exp(-x_i)$ pour $1 \leq i \leq n$, on a que $\mu_{u^0}= - \mu$. Ceci permet de conclure.
\cqfd

\subsubsection{Groupe des automorphismes}\label{Francis}

Rappelons pour commencer que le groupe des automorphismes d'une variété complexe compacte $M$ est un groupe de Lie complexe de dimension finie dont l'algèbre de Lie est $ \eta^\RR(M)$ où $\eta^\RR(M)$ est l'ensemble des champs de vecteurs holomorphes réels de $M$ i.e $Z \in \eta^\RR(M)$ si et seulement si $Z$ est un champ de vecteurs vérifiant $\LLL_Z \omega =0$. De plus, puisque $J$ est intégrable alors $\eta^\RR(M)$ est une algèbre de Lie complexe (voir par exemple \cite{dercalabi}). Si $(M,J)$ une variété complexe compacte, alors on peut voir $T^{1,0}_J M$ comme un fibré vectoriel complexe via l'action de $J$ et on définit $\eta(M)$ comme l'ensemble des champs de vecteurs (complexes) holomorphe i.e. les sections holomorphes du fibré vectoriel $T^{0,1}_JM$. Remarquons que l'application $X \in \eta^\RR(M) \mapsto X^{1,0} \in \eta(M)$ est un isomorphisme d'algèbre de Lie. On pourra consulter \cite{dercalabi} pour plus de détails. \\

Supposons qu'en plus $M$ soit torique de Fano i.e. $(M,\omega,J,g,\TT^n,\mu)$ est une variété kählérienne compacte torique de Fano, et que l'application moment $\mu$ vérifie $\im(\mu)= - P_{alg}$. Cela signifie que l'action de $\TT^n$ s'étend en une action du tore complexe $\TT_\CC$ qui agit de manière fidèle et holomorphe sur $M$ avec une orbite ouverte et dense. En particulier, $\TT_\CC$ définit un tore complexe maximal de $\operatorname{Aut}(M)$. De plus, comme $\TT_\CC$ est connexe, on a que $\TT_\CC$ appartient à  $\operatorname{Aut}^{0}(M)$, donc nous avons 
$$
 \xymatrix{
   ~ \TT^n ~ \ar@{^{(}->}[r]   & ~ \TT_\CC~ \ar@{^{(}->}[r]  & ~ \operatorname{Aut}^{0}(M) ~ \ar@{^{(}->}[r]  &~ \operatorname{Aut}(M)~.  
  }
  $$
  De plus, $\operatorname{Aut}^{\circ}(M)$ admet la décomposition suivante (conséquence directe du corollaire 5.8 de \cite{Fujiki1978}):
$$
\operatorname{Aut}^{\circ}(M) = \operatorname{Aut}(M) \ltimes R_u,
$$
où $\operatorname{Aut}_r(M)$ est un sous-groupe réductif de $\operatorname{Aut}^{\circ}(M)$ et la complexification d'un sous-groupe compacte maximal $K$ de $\operatorname{Aut}^0(M)$ et $R_u$ le radical unipotent de $\operatorname{Aut}^{\circ}(M)$. De plus, si on note $\eta^\RR(M)$, $\eta_r^\RR(M)$ et $\eta_u^\RR(M)$ les algèbres de Lie de $\operatorname{Aut}(M),\operatorname{Aut}_r(M)$ et $R_u$ respectivement, alors
$$
 \xymatrix{
   ~ \ttt ~ \ar@{^{(}->}[r]   & ~ \eta^\RR_0(M)~ \ar@{^{(}->}[r]  & ~\eta^\RR(M) = \eta_r^\RR(M) \oplus \eta_u^\RR(M),  
  }
$$
et via l'isomorphisme $X \mapsto X^{1,0}$, nous obtenons alors une décomposition
$$
 \xymatrix{
   ~ \ttt ~ \ar@{^{(}->}[r]   & ~ \eta_0(M)~ \ar@{^{(}->}[r]  & ~\eta(M) = \eta_r(M) \oplus \eta_u(M),  
  }
$$
Pour décrire $\operatorname{Aut}^0(M)$ et son algèbre de Lie, nous allons introduire la notion de \textit{racines de Demazure}. Notons $\ttt$ l'algèbre de Lie de $\TT^n$.

\bdefi\label{Pierre5}
Un élément $a \in \ttt^*$ est une racine de Demazure de $P_{alg}$ s'il existe un unique $\rho_a \in \Delta(1)$ tel que $\langle a, b_{\rho_a} \rangle=1$ et $\langle a, b_\rho \rangle \leq 0 $ pour tout $\rho \in \Delta(1)$ tel que $\rho \neq \rho_a$. 
\edefi
On note $R(P_{alg})$ l'ensemble des racines de Demazure de $P_{alg}$. On définit alors $$ 
S(P_{alg}) := R(P_{alg}) \cap -R(P_{alg}) = \lbrace a \in R(P_{alg}) ~:~ - a \in R(P_{alg}) \rbrace$$
et $$U(P_{alg})=R(P_{alg}) \backslash S(P_{alg})=\lbrace a \in R(P_{alg}) ~:~ - a \not \in R(P_{alg}) \rbrace.$$ Avec ces notations, nous avons le résultat suivant.
\btheo[Proposition 7 \cite{Demazure1970}]\label{decompoetator}
Nous avons les décompositions suivantes :
$$
\eta(M)= \eta_0(M) \oplus \bigoplus_{\alpha \in R(P_{alg})} \CC \, V_\alpha,
$$
$$
\eta_r(M)= \eta_0(M) \oplus \bigoplus_{\alpha \in S(P_{alg}) } \CC \, V_\alpha,
$$
$$
\eta_u(M)= \bigoplus_{\alpha \in U(P_{alg})} \CC \, V_\alpha,
$$
où pour tout $\alpha=(\alpha_1, \cdots, \alpha_n) \in R(P_{alg})$, le champ de vecteurs $V_\alpha \in \eta(M)$ est un champ de vecteurs défini sur $\TT_\CC= \lbrace (t_1, \cdots,t_n) ~:~ t_i \in \CC^\times \rbrace$ par :
\begin{equation}\label{Valpha}
V_\alpha \vert_{\TT_\CC} = \prod_{i=1}^n t_i^{-\alpha_i} \,  \sum_{k=1}^n (b_{\rho_\alpha})_k \, t_k \dfrac{\partial}{\partial t_k}.
\end{equation}
\etheo
Pour plus de détails, on pourra par exemple consulter la section 3 de \cite{oda2012convex}. 

\section{Solitons de Kähler-Ricci dans le cas torique}

Soit $(X,g)$ un soliton de Kähler-Ricci et $\omega=\omega_g$. Le fait que $X$ soit un champ de vecteurs holomorphe et que $M$ soit une variété de Fano compacte  implique, par la théorie de Hodge, qu'il existe une unique fonction $\theta_X \in \C^{\infty}(M)$
telle que
$$
i_X \omega = \sqrt{-1} \, \overline{\partial} \theta_X, ~~ \int_M e^{\theta_X} \omega_g^n = \int_M \omega_g^n.
$$
Et donc par la formule de Cartan 
$
\LLL_X(\omega)= \sqrt{-1} \, \partial \overline{\partial} \, \theta_X$ (voir par exemple \cite{TZ1,TZ2}). Ainsi, si nous prenons les classes de cohomologie de l'équation \eqref{KRS}, nous avons alors la relation suivante :
\begin{equation*}
\lambda [\omega] = 2 \pi c_1(M).
\end{equation*}
De plus en prenant le produit extérieur de l'équation \eqref{KRS} avec $\omega^{n-1}$ et en intégrant, nous obtenons
\begin{align*}
\int_M \Ric(\omega) \wedge \omega^{n-1} - \lambda \int_M \omega^n &= \sqrt{-1} \int_M \partial \overline{\partial} \theta_X \wedge\omega^{n-1}. \\
\end{align*}
Pour toute $(1,1)$-forme réelle $\alpha$, on définit sa trace par rapport à $\omega$ par la formule :
$
\left( \Tr_{\omega}\alpha  \right) \cdot  \omega^n = 2n \, \alpha \wedge \omega^{n-1},
$
ce qui nous donne
\begin{align*}
\int_M \Tr_\omega(\Ric(\omega)) \omega^n - 2n \, \lambda \, \int_M \omega^n &= \sqrt{-1} \int_M \Tr_\omega ( \partial \overline{\partial} \theta_X ) \omega^{n}.
\end{align*}
Maintenant, par définition, on a que $Scal_g =  \Tr_\omega(\Ric(\omega)) $ et $\Delta^{g}= \Tr_\omega ( \sqrt{-1} \partial \overline{\partial})$, ce qui nous donne
$$
\int_M Scal_g \, \omega^n - 2n \, \lambda \int_M \omega^n = \int_M \Delta^{g} \theta_X \, \omega^n.
$$
En remarquant que le dernier terme est nul par la formule d'intégration par partie, nous obtenons
\begin{equation}
\lambda= \cfrac{1}{2n} \, \overline{Scal},
\end{equation}
où $n$ est la dimension complexe de $M$ et $\overline{Scal} := \int_M Scal_g \, \omega^n / \int_M \omega^n$ la courbure scalaire moyenne de $(M,g)$.

Supposons que $(M,\omega,J,\TT^n, \mu)$  soit une variété kählérienne compacte torique de Fano. On peut montrer (\cite{WZ}) qu'il existe un unique soliton $(X,g)$ de Kähler-Ricci (modulo l'action des automorphismes holomorphes de $M$) et que le champ de vecteurs solitonique est de la forme $$X= JX_a + \sqrt{-1} X_a$$ où $X_a$ est le champ de vecteurs fondamental pour l'action du tore $\TT^n$ associé à un vecteur $a=(a_1,\cdots,a_n) \in \ttt$. En particulier, nous avons 
\begin{equation}\label{thetaXtor}
\LLL_X \omega= -2 \sqrt{-1} \, \partial  \overline{\partial} \, \langle \mu, a \rangle = -\sum_{i,j,k=1}^n a_i H_{ij,k} d x_k \wedge dt_j.
\end{equation}
En effet, on a d'une part
\begin{align*}
\LLL_X \omega & = d i_X \omega
= d [ \omega(JX_a + \sqrt{-1} X_a , \cdot ) ]  
= d [ \omega( J X_a , \cdot )] + \sqrt{-1} \, d [ \omega(X_a, \cdot)] 
= d [ \omega( J X_a , \cdot )] - \sqrt{-1} \, d^2  \langle \mu, a \rangle \\
& = - d [ g( X_a, \cdot ) ] 
 = - d[ \sum_{i,j=1}^n a_i H_{ij} dt_j ]
 = - \sum_{i,j,k=1}^n a_i H_{ij,k} d x_k \wedge dt_j,
\end{align*}
et d'autre part:
\begin{align*}
- 2 \sqrt{-1} \partial \overline{\partial} \langle \mu , a \rangle & = -  d d^c \langle \mu, a \rangle 
= - d ( \sum_{i=1}^n a_i Jd x_i)
= - d (\sum_{i,j=1}^n a_i H_{ij} dt_j) 
= -\sum_{i,j,k=1}^n a_i H_{ij,k} d x_k \wedge dt_j.
\end{align*}
Ainsi, en prenant la trace par rapport à la forme $\omega=\omega_g$ de l'équation \eqref{KRS}, nous sommes ramenés à chercher les couples $(g,a)$ composés d'une métrique kählérienne $g$ sur $M$ et d'un vecteur $a \in \ttt$ vérifiant
\begin{equation}\label{KRSG}
Scal_g - \overline{Scal} = -2 \, \Delta^g \langle \mu, a \rangle.
\end{equation}
Ainsi, sur les variétés toriques, un soliton de Kähler-Ricci peut être vu comme un couple $(g,a)$ vérifiant l'équation \eqref{KRSG}. En fait, on peut montrer (voir \cite{Eve1,Guan1,Guan2}) que cette définition est équivalente dans le cas d'une variété torique kählérienne compacte de Fano à la définition classique d'un soliton.

Terminons en rappelant que le vecteur $a \in \ttt$ est entièrement déterminé par la combinatoire du  polytope $\im(\mu)$ associé à la variété torique symplectique $(M, \omega,\TT^n, \mu)$ grâce à l'annulation de l'invariant de Futaki. Cette dernière condition s'exprime (voir la  fin de la section 2 de \cite{WZ} pour les détails calculatoires) de la façon suivante :
\begin{equation}\label{FutakiTor}
\forall f \in \operatorname{Aff}(\ttt^*,\RR),~~\int_P e^{-2 \langle a,x \rangle} f \, dv(x) =  f(p) \, \int_P e^{-2 \langle a,x \rangle}  \, dv(x),
\end{equation}
où $p$ est le centre privilégié de $P$, $dv$ la forme volume euclidienne standard et $\operatorname{Aff}(\ttt^*,\RR)$ l'ensemble des fonctions affines réelles sur $\ttt^*$. On pourra consulter \cite{Donaldson,WZ}.

\section{Le laplacien pondéré sur une variété torique de Fano }

\subsection{La première valeur propre du laplacien pondéré}\label{s1ervp}

On fixe une variété kählérienne compacte de Fano torique $(M,J,g_0,\omega_0,\TT^n, \mu)$ et un soliton de Kähler-Ricci $(g,a)$. Lorsqu'on travaille sur les solitons de Kähler-Ricci, il est important de pondérer le laplacien comme expliqué dans \cite{Pali2,Pali1}. Dans notre cas, on définit alors \textit{le laplacien pondéré} par la formule suivante :
\begin{equation}\label{lapcomgen}
\forall v \in \C^\infty(M,\RR),~~ \Delta^{g,a} v := \Delta^{g} v + d v \left( \nabla_g  \, \langle \mu ,2a \rangle \right).
\end{equation}
En utilisant les formules \eqref{Laplacien}, \eqref{nablag} et \eqref{gtorique}, on obtient donc l'expression de $\Delta^{g,a}$ dans les coordonnées actions-angles :
\begin{align}\label{LapPond}
\Delta^{g,a}  = - \sum_{i,j=1}^n \cfrac{\partial}{\partial x_i} \left( H_{ij} \cfrac{\partial}{\partial x_j} \right) - \sum_{i,j=1}^n G_{ij} \cfrac{\partial^2}{\partial t_i \partial t_j} + 2 \sum_{i,j=1}^n a_{i} H_{ij} \cfrac{\partial}{\partial x_j}.
\end{align}

De plus, comme nous travaillons sur des variétés kählériennes, il faut aussi tenir compte de la structure complexe $J$ de la variété $M$ en rajoutant un multiple du terme $B_{g,J}^a$ donné par
\begin{equation}\label{Bexp}
\forall v \in \C^\infty(M,\RR),~~ B_{g,J}^a v := g \left(  \nabla_g v , \nabla_\omega \langle \mu, 2 \, a \rangle \right) = d v \left( \nabla_\omega \langle \mu, 2 \, a \rangle  \right).
\end{equation}
\textit{Le laplacien complexe pondéré} est alors donné par la formule 
\begin{equation}\label{lapcom}
\Delta^a_{g,J} := \Delta^{g,a} - \sqrt{-1} B_{g,J}^a.
\end{equation}
De plus, puisque l'opérateur $J$ est antisymétrique pour $g$, on peut étendre la définition précédente par $\CC$-linéarité aux fonctions $u \in \C^\infty(M,\CC)$. En utilisant la formule \eqref{nablaw}, nous obtenons que
\begin{equation}\label{Bexpactang}
B_{g,J}^a = \sum_{j=1}^n 2 \, a_j \cfrac{\partial}{\partial t_j}.
\end{equation}
D'où
\begin{equation}\label{lapcomactang}
\Delta^{a}_{g,J}  = - \sum_{i,j=1}^n \cfrac{\partial}{\partial x_i} \left( H_{ij} \cfrac{\partial}{\partial x_j} \right) - \sum_{i,j=1}^n G_{ij} \cfrac{\partial^2}{\partial t_i \partial t_j} + 2 \sum_{i,j=1}^n a_{i} H_{ij} \cfrac{\partial}{\partial x_j} -  2 \, \sqrt{-1}\sum_{j=1}^n  a_j \cfrac{\partial}{\partial t_j} .
\end{equation}
Terminons en remarque que si $u$ et $v$ sont des fonctions $\TT^n$-invariantes i.e. indépendantes de la variable $t$ alors
\begin{equation}\label{lapuv}
\Delta_{g,J}^a(uv)= v \, \Delta^{g,a}u +  u \, \Delta^{g,a}v - 2 \sum_{i,j=1} \frac{\partial u}{\partial x_i} H_{ij} \frac{\partial v}{\partial x_j}.
\end{equation}

Maintenant, nous nous intéressons aux \textit{valeurs propres} de ce laplacien et aux \textit{fonctions propres} (appelées aussi \textit{potentiels}) associées i.e. on cherche des nombres complexes $\nu \in \CC$ et des fonctions non nulles $f \in \C^{\infty}(M,\CC)$ vérifiant
$$
\Delta^{a}_{g,J} f = \nu \, f .
$$
Le résultat fondamental concernant les valeurs propres du laplacien pondéré complexe $\Delta_{g,J}^a$ est le lemme 28 de \cite{Pali1}. En particulier, il montre que, en renormalisant l'équation \eqref{KRS} pour avoir $\lambda=1$, la première valeur propre est $2$. Pour le compléter, nous avons alors le résultat nouveau suivant qui nous donnent des fonctions propres associées à cette première valeur propre. Il généralise aussi la proposition 2.4 de \cite{Eve2}  qui traite le cas des métriques de Kähler-Einstein sur les variétés toriques.

\blemm\label{1ervp}
Soit $(M,J,g_0,\omega_0,\TT^n, \mu)$ une variété torique kählérienne compacte de Fano. Alors le couple $(g,a)$ est un soliton de Kähler-Ricci si et seulement si à constante  additive près
$$
\forall b \in \ttt,~~2\lambda \langle  \mu, b \rangle = \Delta^{g,a} \langle \mu, b \rangle,
$$
où $\lambda$ est la constante de Kähler-Einstein dans l'équation des solitons \eqref{KRS}. 
\elemm

\dem
La preuve suit la démarche de la proposition 2.4 de \cite{Eve2} en rajoutant les termes nécessaires dans le cas solitonique. 

Soit $(g,a)$ un soliton de Kähler-Ricci. Il existe donc $\lambda \in \RR$ tel que
$$
\Ric(\omega) - \lambda \omega = \LLL_X \omega,
$$
où $\omega$ est la forme de Kähler de $g$ et
$$
X
= 
JX_a + \sqrt{-1}X_a.
$$
Commençons par remarquer, grâce à l'équation \eqref{Laplacien} et \eqref{Ricactang}, que nous avons, pour tout $b \in \ttt$,
\begin{equation}\label{Juju}
d \left[ \Delta^g \langle \mu, b \rangle \right] = d  \left[ \Delta^{g} \left( \sum_{i=1}^n b_i x_i  \right) \right] 
= - d \left[ \sum_{i,j=1}^n b_j H_{ij,i} \right] 
= - \sum_{i,j,k=1}^n b_j H_{ij,ik} dx_k
= -2 i_{X_b}\Ric(\omega).
\end{equation}
De plus, nous avons aussi (voir les formules \eqref{thetaXtor} et \eqref{nablag}) :
\begin{align*}
i_{X_b}( \LLL_X \omega  )
=   \sum_{i,j,k=1}^n a_i H_{ij,k} b_j dx_k 
= d \left[d \langle \mu , b \rangle \left( \nabla_g \langle \mu, a \rangle \right) \right].
\end{align*}
Ainsi nous obtenons grâce à l'équation des solitons \eqref{KRS} et les formules \eqref{fanfan1}, \eqref{Juju}
\begin{align*}
d \Delta^{g,a} \langle \mu , b \rangle = i_{X_b} \left( -2 \Ric(\omega) + 2 \LLL_X \omega \right)
= -2  \lambda i_{X_b} \omega 
= 2 \lambda d \langle \mu,b \rangle 
\end{align*}
On conclut alors en utilisant le fait que la variété $M$ est compacte connexe qu'il existe une constante $c_b \in \CC$ telle que
$$
\Delta^{g,a} \langle \mu, b \rangle =  2 \lambda ( \langle \mu, b \rangle  + c_b ),
$$
ainsi nous obtenons bien
$$
\forall b \in \ttt,~~ \exists c_b \in \CC,~~\Delta^{g,a} \left( \langle \mu, b \rangle + c_b \right) =  2 \lambda ( \langle \mu, b \rangle  + c_b).
$$

Réciproquement, supposons cette propriété vérifiée, on a, en utilisant la formule \eqref{lapcomgen}, pour tout $i=1 ,\cdots,n $, où $b(i)$ a toutes ses coordonnées nulles sauf la $i$-ème,
\begin{equation}\label{tildeHa}
2 \lambda \, (x_i + c_{b(i)}) = -\sum_{j=1}^n H_{ji,j} + 2 \sum_{j=1}^n a_j H_{ji}.
\end{equation}
Ce qui nous donne, par les formules \eqref{Ricactang} et \eqref{Laplacien},
\begin{align*}
\Ric(\omega)&= -\cfrac{1}{2} \sum_{i,k,l=1}^n  H_{il,ik} \, d x_k \wedge d t_l = - \cfrac{1}{2} \sum_{i,k,l=1}  \cfrac{ \partial  H_{il,i}}{\partial x_k} \, d x_k \wedge d t_l \\
&= \cfrac{1}{2} \sum_{k,l=1}  \cfrac{ \partial}{\partial x_k} \left[ 2 \lambda \, (x_l + c_{b(l)}) -   2 \sum_{i=1}^n a_i H_{il } \right] d x_k \wedge d t_l \\
&= \lambda \sum_{k=1}^n d x_k \wedge dt_k -  \sum_{l,k,i=1}^n a_i H_{il,k} \, d x_k \wedge dt_l = \lambda \omega + \LLL_X \omega.
\end{align*}
\cqfd
\medskip

On peut étendre le lemme \ref{1ervp} au cas complexe grâce au corollaire suivant :

\bcoro\label{proprecst=0}
Pour tout $b \in \ttt$, la fonction $ \langle \mu,b \rangle$ (modulo une constante additive dépendant de $b$) est encore une fonction propre de valeur propre $2\lambda$ pour le laplacien pondéré complexe $\Delta^a_{g,J}$.
\ecoro
\dem
Il suffit de remarquer que $B^a_{g,J} \langle \mu, b\rangle=0$ pour tout $b \in \ttt$ par l'équation \eqref{Bexpactang} car la fonction $ \langle \mu, b \rangle$ ne dépend pas des coordonnées $t_1, \cdots, t_n$.
\cqfd
\medskip

Lorsque le polytope est $-P_{alg}$ alors on peut se débarasser de la constante grâce au résultat suivant.

\bcoro\label{corro1er}
Soient $(M,\omega_0,J,\TT^n, \mu)$ une variété torique kählérienne compacte de Fano telle que $\im(\mu)= -P_{alg}$ et $(g,a)$ un soliton de Kähler-Ricci. Alors pour tout $b \in \ttt$, la fonction $\langle \mu, b \rangle$ est une fonction propre de valeur propre $2 \lambda$ pour le laplacien pondéré complexe $\Delta^a_{g,J}$.
\ecoro

\dem
Pour tout $b \in \ttt$, on sait, par le corollaire \ref{corro1er} qu'il existe une constante $c_b$ telle que
$$
\Delta_{g,J}^a \langle \mu , b \rangle = 2 \lambda \, \langle \mu , b \rangle +2 c_b.
$$
On intègre alors par rapport à $e^{\theta_X} \omega_g^n$ où $X= JX_a + \sqrt{-1} \, X_a$ et on obtient en faisant une intégration par parties que
$$
0= \int_M \Delta_{g,J}^a \langle \mu , b \rangle \, e^{\theta_X} \omega_g^n  =  2 \, \lambda \, \int_M \langle \mu , b \rangle \, e^{\theta_X} \omega_g^n + 2 \, c_b \int_M e^{\theta_X} \omega_g^n.
$$
Pour conclure, il suffit donc de montrer que
$$
\int_M \langle \mu, b \rangle \, e^{\theta_X} \omega_g^n  = 0.
$$
Puisque que $\theta_X = -2 \langle \mu, a \rangle + \tilde{c}_a$ (voir l'équation \eqref{thetaXtor}) où $\tilde{c}_a$ est une constante, il suffit de montrer que
$$
\int_M \langle \mu, b \rangle \, e^{-2 \langle \mu, a \rangle} \omega_g^n  = 0.
$$
En notant $f$ le potentiel kählérien de $\omega_g=2 \sqrt{-1} \, \partial \overline{\partial} \, f $, par l'équation \eqref{Paulin1}, cela revient à montrer que
$$
\int_{\RR^n}  \sum_{k=1}^b b_k \, \dfrac{\partial f}{\partial u_k}  \, \exp \left( {-2 \sum_{l=1}^n a_l \dfrac{\partial f}{\partial u_l}} \right) \, \det(f_{ij}) \, du = 0
$$
Et en faisant le changement de variable $x_i = \dfrac{\partial f}{\partial u_i}$, nous devons finalement montrer que
$$
\int_{\im(\mu)}  \sum_{k=1}^b b_k \, x_k \, \exp \left( -2 \sum_{l=1}^n a_l x_l \right) \, dx = 0.
$$
Cette dernière intégrale est bien nulle car elle correspond à l'annulation de l'invariant de Futaki (voir l'équation \eqref{FutakiTor}) en remarquant que le centre privilégié de $\im(\mu)=-P_{alg}$ est l’origine par définition de $P_{alg}$.
\cqfd

\section{Étude de la décomposition solitonique}

On fixe une variété torique kählérienne compacte de Fano $(M,J,g_0,\TT^n,\mu)$ et on note $(X,g)$ l'unique soliton de Kähler-Ricci dont on note  $a \in \ttt$ l'élément tel que $X= JX_a + \sqrt{-1} \, X_a$. De plus, on suppose que la constant d'Einstein est normalisé pour avoir $\lambda=1$.

\subsection{Rappel sur la décomposition solitonique}

Avant d'énoncer le théorème, rappelons que si $E$ est un espace vectoriel complexe alors on note $\overline{E}$ le conjugué complexe de cet espace vectoriel complexe. De plus, on pondère le crochet de Poisson usuel de la façon suivante :
$$
\lbrace u, v \rbrace^X_\omega := \lbrace u, v \rbrace_\omega - \int_M \lbrace u, v \rbrace_\omega e^{\theta_X} \omega_g^n.
$$
Le théorème originel est issue de \cite{Pali1,dercalabi} mais nous le donnons ici directement avec les notations des variétés toriques.
\btheo\label{decomposolito}
On a
\begin{enumerate}
\item[$\bullet$] L'application suivante entre algèbres de Lie 
$$
\begin{array}{lrcl}
\chi : & \left( \, \overline{\ker(\Delta^a_{a,J} - 2 \II)}, \sqrt{-1} \, \lbrace \cdot, \cdot\rbrace^X_\omega \, \right)&
\stackrel{\sim}{\longrightarrow} & \left(\eta^\RR(M), [\cdot,\cdot] \right) \\
   & u & \longmapsto  &  \nabla_{g,J} u := \nabla_g \Re(u) + J \, \nabla_g \Im(u). \end{array}
$$
est bien définie et est un isomorphisme d'algèbres de Lie. En particulier, $\ker(\Delta^a_{g,J} - 2 \II)$ est de dimension finie.
\item[$\bullet$] La première valeur propre $\lambda_1(\Delta^a_{g,J})$ de l'opérateur $\Delta^a_{g,J}$ vérifie $\lambda_1(\Delta^a_{g,J}) \geq 2$ avec égalité si $\eta^\RR(M) \neq 0$.
\item[$\bullet$] L'application suivante est bien définie et définit un isomorphisme d'espaces vectoriels réels :
$$
 J \nabla_g : \ker(\Delta^a_{g,J}-2 \II) \cap \C^\infty_X(M,\RR)_0\longrightarrow \kappa(M),
$$
où $\kappa(M)$ est l'ensemble des champs de vecteurs de Killing de la variété riemanienne $(M,g)$ et $$\C^\infty_X(M,\RR)_0 = \lbrace f \in \C^\infty(M,\RR) ~:~ \int_M f e^{\theta_X} \omega_g^n=0 \rbrace.$$
\item[$\bullet$] La forme hermitienne définie sur $\overline{\ker(\Delta^a_{g,J} - 2 \II)}$ par
$$
(u,v) \mapsto \int_M \sqrt{-1} \lbrace u , \overline{v} \rbrace_\omega^X e^{\theta_X} \omega^n_g,
$$
est positive. Et si on note $(\gamma_i)_{i=0,\cdots,N}$ son spectre alors on a la décomposition suivante :
$$
\eta^\RR(M):= \bigoplus_{i=0}^N V_{\gamma_i},
$$
où 
$$
V_{\gamma_i}:= \lbrace \xi \in \eta^\RR(M) ~:~ [\xi, \nabla_g \theta_X] = \gamma_i \xi \rbrace .
$$
En particulier, on a $\gamma_0=0$ et 
$$
V_0 = \kappa(M) + J \, \kappa(M).
$$
\end{enumerate}
\etheo

\subsection{Le sous-espace $\eta_0^\RR$}

On veut étudier l'isomorphisme $\chi$ donné par le théorème \ref{decomposolito}. Rappelons que le tore compact $\TT^n$ s'étend en un tore complexe $\TT_\CC$ qui s'injecte dans le groupe $\operatorname{Aut}^0(M)$, on note $\iota : \TT_\CC \rightarrow \operatorname{Aut}^0(M)$ cette injection. Ainsi, si $\ttt$ est l'algèbre de Lie de $\TT$ et $\ttt_\CC$ l'algèbre de Lie de $\TT_C$, nous obtenons que $\eta_0^\RR(M) = d\iota(\ttt_\CC)$ est une sous-algèbre de Lie de $\eta^\RR(M)$ ( où $\eta^\RR(M)$ est définie dans la section \ref{Francis}) et qu'elle est aussi égale à la complexification de $d \iota(\ttt)$ i.e. $\eta_0^\RR(M) = d \iota(\ttt) \oplus J \, d \iota(\ttt)$. En particulier, on a que $\eta_0^\RR(M)$ est isomorphe à $\ttt \oplus \sqrt{-1} _, \ttt$.

Nous avons alors le lemme suivant qui caractérise la restriction de $\chi^{-1}$ à $\eta_0^\RR(M)$.

\blemm\label{eta0decompo}
L'application $\chi^{-1}$ restreinte à $ \eta_0^\RR(M)$ est donnée via l'isomorphisme $\iota$  par
$$
\begin{array}{l|rcl}
\chi^{-1} \vert_{\eta_0^\RR} : & \eta_0^\RR  & \longrightarrow & \operatorname{Aff}_0^\CC(P) \\
  & d \iota(b_1) + J \, d \iota(b_2) & \longmapsto  & - \langle \mu, b_2 \rangle + \sqrt{-1} \langle \mu, b_1 \rangle \\
   \end{array},
$$
où $\operatorname{Aff}_0^\CC(P)$ est l'ensemble des fonctions lisses  définies sur $M$ s'écrivant sur $M^0$ et en coordonnées actions-angles sous la forme $ \langle \mu, b_1 \rangle + \sqrt{-1} \langle \mu, b_2 \rangle$ pour $(b_1,b_2) \in \ttt \times \ttt$.
\elemm

\dem Commençons par remarquer que $\operatorname{Aff}_0^\CC(P)$ est bien inclus dans $\ker(\Delta^a_{g,J} - 2 \II)$ grâce au corollaire \ref{proprecst=0}. Grâce à l'action du tore, $ b_1 \in \ttt$ induit un champ de vecteurs $X_{b_1}$ qui se trouve être égal à $ \nabla_\omega  \langle \mu, b_1 \rangle $ (voir formule \eqref{nablaw}). Puisque $g=\omega(\cdot, J \cdot)$, ce dernier est égal à $ J \, \nabla_g \langle \mu, b_1 \rangle$. Ainsi si $b_2 \in \ttt$ alors $ \sqrt{-1} \, b_2$ induit le champ de vecteurs $ - \nabla_g \langle \mu, b_2 \rangle$. Ceci permet de conclure en utilisant l'isomorphisme $\chi$ donné par le théorème \ref{decomposolito}.
\cqfd
\medskip

Pour le reste de l'espace $\eta^\RR(M)$, il est préférable de travailler avec les champs de vecteurs holomorphes complexes, nous allons donc transporter la décomposition solitonique dans cet espace et l'étudier ensuite.

\subsection{La décomposition solitonique complexe}
Commençons par remarquer que nous pouvons traduire l'isomorphisme pour les champs de vecteurs holomorphes. En effet, l'algèbre de Lie $\eta^\RR(M)$ est isomorphe à $\eta(M)$ via l'application $X \mapsto X^{1,0}$. Nous obtenons donc l’isomorphisme
$$
\begin{array}{l|rcl}
\tilde{\chi} : & \overline{\ker(\Delta^a_{g,J}- 2 \II )} & \longrightarrow & \eta(M) \\
  & v & \longmapsto  &  (\nabla_{g,J}v)^{1,0}. \\
   \end{array}
$$
De plus, d'après le théorème \ref{decompoetator}, nous avons
$$
\eta(M)= \eta_0(M) \oplus \bigoplus_{\alpha \in R(P_{alg})} \CC \, V_\alpha.
$$
Or, l'espace $\eta_0(M)$, qui est isomorphe à $\eta^\RR_0(M)$ via $X \mapsto X^{1,0}$ étant déjà compris, grâce au lemme \ref{eta0decompo}, on est donc ramené à rechercher une fonction $v_\alpha \in \C^\infty(M,\CC)$ telle que
\begin{equation}\label{systdecompotor}
\left \{
\begin{array}{c @{=} c}
    (\nabla_{g,J}v_\alpha)^{1,0} & V_\alpha \\
    (\Delta^a_{g,J} \overline{v_\alpha} - 2 \, \overline{v_\alpha} )  & 0. \\
\end{array}
\right.
\end{equation}

\subsubsection{Résolution locale du système \eqref{systdecompotor}}

Maintenant, on considère les coordonnées logarithmiques sur $\TT_\CC$ i.e. $t_i= \exp(z_i)$ où $z_i=u_i + \sqrt{-1} v_i$. Dans ce système, nous obtenons, par l'équation \eqref{Valpha}, que
$$
V_\alpha \vert_{\TT_\CC} = \exp \left(\langle -\alpha, z \rangle \right) \, \sum_{k=1}^n (b_{\rho_\alpha})_k \,\dfrac{\partial}{\partial z_k},
$$
où la notation $\langle \alpha, z \rangle$ désigne, avec $ \alpha = ( \alpha_1, \cdots, \alpha_n)$,
$$
\langle \alpha, z \rangle = \sum_{i=1}^n \alpha_i \left( u_i + \sqrt{-1} v_i \right).
$$
On cherche donc une fonction $v_\alpha \in \C^\infty(M,\CC)$ solution du système d'équations \eqref{systdecompotor}. Commençons par remarquer que nous avons
$$
\omega( (\nabla_{g,J}v_\alpha)^{1,0}, \cdot )  = \sqrt{-1} \, \overline{\partial} v_\alpha,
$$
et donc nous devons avoir
$$
\omega(V_\alpha, \cdot ) = \sqrt{-1} \, \overline{\partial} v_\alpha.
$$
Si on se place dans le système de coordonnées actions-angles, nous obtenons donc que
$$
\forall j=1,\cdots,n,~~\dfrac{\partial v_\alpha}{\partial \overline{z}_j} = 2 \sum_{i=1}^n \dfrac{\partial^2 f}{\partial z_i \partial \overline{z}_j} \left( \exp \left(\langle -\alpha, z \rangle \right) \, (b_{\rho_\alpha})_i  \right),
$$
où on rappelle que $f$ est le potentiel kählérien de $\omega$ (voir l'équation \eqref{potf}).
En particulier, cela nous dit que
\begin{align*}
v_\alpha \vert_{M^0}& = 2 \sum_{i=1}^n\dfrac{\partial f}{ \partial z_i} \, (b_{\rho_{\alpha}})_i \,  \exp \left[\langle -\alpha, z \rangle \right] +v_\alpha^0  \\
&= \sum_{i=1}^n \dfrac{\partial f}{ \partial u_i} \, \, (b_{\rho_{\alpha}})_i \,  \exp \left[\langle -\alpha, z \rangle \right] +v_\alpha^0 &   \text{ (puisque le potentiel kählérien est $\TT^n$-invariant)}\\
&= \langle \mu, b_{\rho_{\alpha}} \rangle \, \exp\left(\langle -\alpha, \nabla \phi + \sqrt{-1} t \rangle \right) +v_\alpha^0 & \text{ (par l'équation \eqref{Pierre3})}
\end{align*}
où $\phi$ est le potentiel symplectique (voir le théorème \ref{potsymp}) de $g$ et $v_\alpha^0$ est une fonction holomorphe sur $M^0$.

Il reste à déterminer $v_\alpha^0$. Pour cela, on utilise la seconde équation de \eqref{systdecompotor}. On pose alors
$$
v'_\alpha = \langle \mu, b_{\rho_{\alpha}} \rangle \, e^{\langle - \alpha , \nabla \phi + \sqrt{-1} t \rangle}
$$
et on calcule
$$
\Delta_{g,J}^a \overline{v'_\alpha} - 2 \overline{v'_\alpha}.
$$
Remarquons que nous avons, grâce à l'équation \eqref{lapuv} et au corollaire \ref{corro1er},
\begin{align}
\Delta_{g,J}^a \overline{v'_\alpha} &= \Delta_{g,J}^a \left(\langle \mu, b_{\rho_{\alpha}} \rangle \, \exp\langle -\alpha, \nabla \phi - \sqrt{-1} t \rangle \right) \nonumber \\
& = \Delta_{g,J}^a \left(\langle \mu , b_{\rho_{\alpha}} \rangle \, \exp \langle -\alpha, \nabla \phi  \rangle  \right)  \exp \langle \alpha, \sqrt{-1} t \rangle 
+ \left( \langle \mu, b_{\rho_{\alpha}} \rangle \, \exp \langle -\alpha, \nabla \phi  \rangle  \right) \Delta_{g,J}^a \left(  \exp \langle \alpha,\sqrt{-1} \, t \rangle  \right) \nonumber \\
&= \Delta_{g,J}^a \left(\langle \mu , b_{\rho_{\alpha}} \rangle ) \, \exp \langle -\alpha, \nabla \phi  \rangle  \right)  \exp \langle \alpha, \sqrt{-1} t \rangle 
+
 \left(\langle \mu , b_{\rho_{\alpha}} \rangle \right) \, \Delta_{g,J}^a \left( \exp \langle -\alpha, \nabla \phi  \rangle  \right)  \exp \langle \alpha, \sqrt{-1} t \rangle 
&  \nonumber\\ 
&~~~~+ 
\left( \langle \mu, b_{\rho_{\alpha}} \rangle \, \exp \langle -\alpha, \nabla \phi  \rangle  \right)  \Delta_{g,J}^a \left(  \exp \langle \alpha,\sqrt{-1} \, t \rangle  \right) \nonumber \\
&= 2 \, \left(\langle \mu , b_{\rho_{\alpha}} \rangle ) \, \exp \langle -\alpha, \nabla \phi  \rangle  \right)\exp \langle \alpha, \sqrt{-1} t \rangle 
+
 \left(\langle \mu , b_{\rho_{\alpha}} \rangle \right) \, \Delta_{g,J}^a \left( \exp \langle -\alpha, \nabla \phi  \rangle  \right) \exp \langle \alpha, \sqrt{-1} t \rangle 
& \nonumber\\ 
&~~~~+ 
\left( \langle \mu, b_{\rho_{\alpha}} \rangle \, \exp \langle -\alpha, \nabla \phi  \rangle  \right)  \Delta_{g,J}^a \left(  \exp \langle \alpha,\sqrt{-1} \, t \rangle  \right). \label{Seb1}
\end{align}
De plus, nous avons, par l'équation \eqref{lapcomactang} et puisque $\langle \alpha + \sqrt{-1} t \rangle$ ne dépend pas de $x_1,\cdots,x_n$,
\begin{equation}\label{Seb3}
\Delta_{g,J}^a \left(  \exp\langle \alpha,\sqrt{-1} \, t \rangle  \right) = \left( ^t \alpha G \alpha + 2 \langle \alpha , a\rangle \right) \exp  \langle \alpha, \sqrt{-1} t \rangle.
\end{equation}
En utilisant à nouveau l'équation \eqref{lapcomactang} et en se souvenant que la matrice de $G$ est la matrice hessienne de $\phi$,
\begin{align}\label{Seb2}
\Delta_{g,J}^a \left( \exp \langle -\alpha, \nabla \phi  \rangle \right) = - \left( ^t \alpha G \alpha + 2 \langle \alpha , a\rangle \right) \, \exp \langle -\alpha, \nabla \phi  \rangle. 
\end{align}
Nous avons donc finalement, à l'équation \eqref{Seb1}, que
\begin{equation}
\Delta_{g,J}^a \overline{v'_\alpha} = 2 \overline{v'_\alpha} + 2 \langle \alpha, b_{\rho_\alpha} \rangle  \exp \langle -\alpha, \nabla \phi - \sqrt{-1} t  \rangle.
\end{equation}
Or $\langle \alpha, b_{\rho_\alpha} \rangle =1$ (puisque $\alpha$ est une racine de Demazure, voir la définition \ref{Pierre5}) donc nous avons finalement :
$$
\Delta_{g,J}^a \overline{v'_\alpha} - 2 \overline{v'_\alpha} = 2 \exp \langle -\alpha, \nabla \phi - \sqrt{-1} t  \rangle.
$$
Remarquons alors que, grâce aux équations \eqref{Seb3} et \eqref{Seb2},
$$
\Delta_{g,J}^a \left( \exp \langle -\alpha, \nabla \phi - \sqrt{-1} t  \rangle \right) =0.
$$
Donc finalement
\begin{align*}
\Delta_{g,J}^a \left( \overline{v'_\alpha} +\exp \langle -\alpha, \nabla \phi - \sqrt{-1} t  \rangle \right) &=  2 \left( \overline{v'_\alpha} +  \exp \langle -\alpha, \nabla \phi - \sqrt{-1} t \rangle \right).
\end{align*}

Nous posons alors la fonction suivante :
$$
v_\alpha := \left( \langle \mu, b_{\rho_{\alpha}} \rangle +1 \right) \, \exp\langle -\alpha, \nabla \phi + \sqrt{-1} t \rangle \in \C^{\infty}(M^0,\CC).
$$
On l'appellera \textit{la fonction de racine $\alpha$}.
\btheo
La fonction de racine $\alpha$ est une solution locale i.e. sur $M^0$ au système d'équations \eqref{systdecompotor}.
\etheo
\dem
La seconde équation de \eqref{systdecompotor} est vraie d'après ce qui précède, il s'agit donc de regarder si le terme additif
$$
v_\alpha^0=\exp\langle -\alpha, \nabla \phi + \sqrt{-1} t \rangle 
$$
n'interfère pas dans la résolution de la première équation. Or en coordonnées logarithmiques ce terme s'écrit 
$$
v_\alpha^0=\exp\langle -\alpha, z \rangle. 
$$
 Cette fonction est clairement holomorphe sur $M^0$. Ceci permet de conclure.
\cqfd

\subsubsection{Résolution globale du système \eqref{systdecompotor}}

On veut maintenant trouver une fonction appartenant à $\C^\infty(M,\CC)$ qui vérifie le système d'équations \eqref{systdecompotor}. On va donc montrer le résultat suivant.
\btheo
La fonction $v_\alpha \in \C^\infty(M^0,\CC)$ s'étend en une fonction appartenant à $\C^\infty(M,\CC)$.
\etheo
\dem
Dans un premier temps, nous allons montrer que $v_\alpha$ est une fonction continue sur la variété $M$. On sait que la fonction $v_\alpha$ a pour expression
$$
v_\alpha = \left( \langle \mu, b_{\rho_{\alpha}} \rangle +1 \right) \, \exp\langle -\alpha, \nabla \phi + \sqrt{-1} t \rangle. 
$$
Or la théorie des potentiels symplectiques (voir le théorème \ref{potsymp}), nous dit que $\phi$ peut s'écrire
$$
\phi = \phi_0 + h,
$$
où $\phi_0$ est le potentiel de Guillemin (voir l'équation \eqref{Guillemin}) et $h \in \C^\infty(U,\RR)$ avec $U$ ouvert de $\RR^n$ contenant $P$. Ainsi, on peut écrire
$$
v_\alpha = \left( \langle \mu, b_{\rho_{\alpha}} \rangle +1 \right) \, \exp\langle -\alpha, \nabla \phi_0 \rangle \cdot   \exp \langle -\alpha, \nabla h + \sqrt{-1} t \rangle.
$$
Remarquons tout de suite que $\exp \langle -\alpha, \nabla h + \sqrt{-1} t \rangle$ s'étend en une fonction continue sur $P \times \TT^n$. Pour le terme restant, en utilisant l'expression du potentiel de Guillemin (voir l'équation \eqref{Guillemin}), on obtient que
\begin{align*}
\left( \langle \mu, b_{\rho_{\alpha}} \rangle +1 \right) \, \exp\langle -\alpha, \nabla \phi_0 \rangle &=  e^{ - \frac{1}{2}\langle \alpha, \sum_{\rho \in \Delta(1)} b_\rho \rangle}\left( \langle \mu, b_{\rho_{\alpha}} \rangle +1 \right)  \, \prod_{\rho \in \Delta(1)} \left( \langle \mu , b_\rho \rangle + 1 \right)^{- \frac{1}{2}\langle \alpha, b_\rho \rangle}.
\end{align*}
Par les propriétés des racines de Demazure (voir la définition \ref{Pierre5}), nous avons
$$
\langle \alpha, b_{\rho_{\alpha}} \rangle = 1,~~\langle \alpha, b_{\rho}  \rangle \leq 0 ~~ \forall \rho \neq \rho_\alpha.
$$
Ceci permet d'étendre (par la même formule) de manière continue la fonction $\left( \langle \mu, b_{\rho_{\alpha}} \rangle +1 \right) \, \exp\langle -\alpha, \nabla \phi_0 \rangle$ au bord du polytope. Nous avons donc montré que $v_\alpha$ définit une fonction continue sur $P \times \TT^n$. De plus, nous avons
$$
\forall t \in \TT^n,~~\forall F \in \F(P),~~ \forall x \in F,~~ \forall \rho \in I_F,~~ v_\alpha(x, b_\rho +t) = v_\alpha(x,t).
$$
En effet, il y a deux cas à traiter :
\begin{itemize}
\item[$\bullet$] Supposons que $\langle \alpha, b_\rho \rangle =0$ alors par définition la fonction $v_\alpha$ vérifie la propriété demandée.
\item[$\bullet$] Il reste à traiter le cas où $\langle \alpha, b_\rho \rangle  \neq 0$. Dans ce cas, la fonction $v_\alpha$ est nulle sur la face correspondante donc la propriété est vérifiée.
\end{itemize}
Ainsi, nous obtenons que $v_\alpha$ induit (par la même expression) une fonction continue sur $F \times \TT^n/\TT_F$ et donc nous obtenons que $v_\alpha$ définit une fonction continue sur l'espace $\bigsqcup_{F \in \F(P)} F \times \TT^n/\TT_F$.

Pour conclure, on sait (voir par exemple la section 2.1 de \cite{Eve3}) que l'espace topologique sous-jacent à la variété $M$ peut se décrire comme
$$
M \simeq \left( \bigsqcup_{F \in \F(P)} F \times \TT^n/\TT_F \right) ~/~ \sim
$$
où $(x,t) \in F \times \TT^n/\TT_F \sim (x',t') \in F' \times \TT^n/\TT_{F'}$ si et seulement si
 $x=x'$
et
$ t = t' \mod \TT_{F \cap F'}$. Il faut et suffit donc de montrer que si $(x,t) \in F \times \TT^n/\TT_F \sim (x',t') \in F' \times \TT^n/\TT_{F'}$  alors $v_\alpha(x,t)=v_\alpha(x,t')$. Ce qui est le cas par définition de la fonction $v_\alpha$.

Maintenant, on va montrer que la fonction $v_\alpha$ est lisse. Pour cela, on remarque qu'elle définit une fonction lisse sur $M^0$ qui vérifie 
$$
(\nabla_{g,J}v_\alpha)^{1,0} = V_\alpha \text{ sur $M^0$.}
$$
Or on sait (par la théorie de Hodge) qu'il existe une fonction $\theta_{V_\alpha} \in \C^{\infty}(M,\CC)$ telle que
$$
(\nabla_{g,J} \theta_{V_\alpha})^{1,0} = V_\alpha \text{ sur $M$},
$$
et donc par restriction sur $M^0$. Cela signifie que sur $M^0$, nous avons 
$$
\overline{\partial} \left[\theta_{V_\alpha} - v_\alpha \right] = 0.
$$
Cela signifie qu'il existe une fonction holomorphe $\theta_\alpha^0$ définie sur $M^0$ tel que
$$
\theta_{V_\alpha} - v_\alpha = \theta_\alpha^0.
$$
Le membre de gauche étant borné puisque continu sur la variété compacte $M$, nous obtenons que $\theta_\alpha^0$ est une constante et donc par densité nous obtenons que
$
\theta_{V_\alpha} = v_\alpha
$
modulo une constante additive sur $M$ et donc $v_\alpha$ est aussi lisse.
\cqfd
\medskip

Nous pouvons résumer ce que nous venons de démontrer dans le théorème suivant.
\btheo
Soit $(M,\omega_0,J,\TT^n,\mu)$ une variété torique kählérienne compacte de Fano telle que $P:=\im(\mu)$ soit le polytope algébrique associé dont on note $R(P)$ l'ensemble de ses racines de Demazure. Supposons que $(g,a)$ soit un soliton de Kähler-Ricci tel que sa constante de Kähler-Einstein soit égale à $1$ sur la variété complexe $(M,J)$. 
\begin{itemize}
\item[$\bullet$] Nous avons la décomposition suivante :
$$
\overline{\ker \left( \Delta^a_{g,J} - 2 \II \right)} = \operatorname{Aff}_0^\CC \oplus \bigoplus_{\alpha \in R(P)} \CC \, \tilde{v_\alpha},
$$
où $\operatorname{Aff}_0^\CC$ est l'ensemble des fonctions lisses sur $M$ s'écrivant sous la forme $\langle x , b_1 \rangle + \sqrt{-1} \langle x , b_2 \rangle$ où $(b_1,b_2) \in \ttt \times \ttt$ et $\tilde{v_\alpha}$ est la fonction appartenant à $\C^\infty(M,\CC)$ telle qu'en coordonnées action-angles (x,t)
$$
\tilde{v_\alpha} \vert_{M^0} = \left( \langle x, b_{\rho_\alpha} \rangle + 1 \right) e^{- \langle \alpha , \nabla \phi - \sqrt{-1} \, t \rangle },
$$
où $\phi$ est le potentiel symplectique associé à $g$.

\item[$\bullet$] De plus, si on note $V_{\gamma_i}$ les espaces propres de la décomposition solitonique et $\chi$ l'isomorphisme entre $ \overline{\ker \left( \Delta^a_{g,J} - 2 \II \right)}$ et $\eta^{\RR}(M)$ du théorème \ref{decomposolito}, alors nous avons
$$
\chi^{-1}(V_0)= \operatorname{Aff}_0^\CC \oplus \bigoplus_{\alpha \in R(P), \langle \alpha, a \rangle=0} \CC \, \tilde{v_\alpha}
$$
et pour tout $i \geq 1$,
$$
\chi^{-1}(V_{\gamma_i})= \bigoplus_{\alpha \in R(P), \, 2 \langle \alpha, a \rangle=\gamma_i} \CC \, \tilde{v_\alpha}.
$$
\item[$\bullet$] En particulier, les fonctions $\tilde{v_\alpha}$ sont des fonctions propres de l'opérateur $\Delta^a_{g,-J} - 2 \II $ pour la valeur propre $4 \langle \alpha , a \rangle$.
\end{itemize}
\etheo
\dem
Le premier point est une conséquence de la discussion qui précède. 
Le second point est une conséquence de l'isomorphisme $\chi$ énoncé dans le théorème \ref{decomposolito} et du fait que les $\overline{\tilde{v_\alpha}}$ sont bien des solutions au système d'équations \eqref{systdecompotor}.
Le troisième point est alors un calcul direct.
\cqfd

\section{Exemples de décompositions solitoniques dans le cas torique}

\subsection{L'espace projectif}

L'exemple le plus simple est de considérer l'espace projectif. Pour encore plus de simplicité, on va se limiter à $\CC \PP^2$ avec sa structure complexe standard, sa métrique de Fubini-Study et l'action de $\TT^n$ donnée par $(t_1,t_2) \cdot [z_0,z_1, z_2] = [z_0, t_1 \, z_1, t_2 \, z_2]$. Le polytope de Delzant est alors donné (voir \cite{Abreu}) par 
$$
P = \bigcap_{i=1}^3 \left\lbrace x \in \RR^2 ~:~ \langle x , b_i \rangle +1 \geq 0 \right\rbrace
$$
où $\RR^2$ est muni de sa base canonique $(e_1,e_2)$ et de son produit scalaire usuel permettant de l'identifier avec son dual, et
$$b_1:=e_1~~,b_2:=e_2,~~ b_3:=-e_1 -e_2.$$

De plus, on sait que la métrique de Fubini-Study de $\CC \PP^2$ est une métrique de Kähler-Einstein dont le potentiel symplectique est le potentiel de Guillemin $\phi_0$ donné par
$$
\phi_0 = \cfrac{1}{2} \, \sum_{r=1}^3 l_r \ln l_r,
$$
où on a posé $l_r := \langle x , b_r \rangle +1$. On pourra consulter \cite{Eve2} ou \cite{Abreu}. Nous obtenons alors que
$$
\nabla \phi_0 =\cfrac{1}{2} \, \sum_{r=1}^3 \left( 1 + \ln l_r \right) b_r = \cfrac{1}{2} \, \ln \frac{x_1+1}{1-x_1-x_2} \, e_1 + \cfrac{1}{2} \, \ln \frac{x_2+1}{1-x_1-x_2} \, e_2.
$$
On veut maintenant déterminer la décomposition solitonique. Il faut donc déterminer les racines de Demazure :
\begin{enumerate}
\item On cherche $\alpha \in \ZZ^2$ tel que
\[
\left \{
\begin{array}{cc}
    \langle \alpha, e_1 \rangle & = 1 \\
      \langle \alpha, e_2 \rangle & \leq 0 \\
        \langle \alpha, -e_1-e_2 \rangle & \leq 0 \\
\end{array}
\right.
\]
On obtient deux racines :
$$
\alpha_{1,1}:= e_1, ~~ \alpha_{1,2}:=e_1-e_2.
$$
\item De même on obtient :
$$
\alpha_{2,1}:= e_2, ~~ \alpha_{2,2}:= -e_1 + e_2.
$$
\item On doit terminer en résolvant :
\[
\left \{
\begin{array}{cc}
    \langle \alpha, e_1 \rangle & \leq 0 \\
      \langle \alpha, e_2 \rangle & \leq 0 \\
        \langle \alpha, -e_1-e_2 \rangle & = 1 \\
\end{array}
\right.
\]
On obtient aussi deux racines :
$$
\alpha_{3,1}:=-e_1, ~~ \alpha_{3,2}:= -e_2.
$$
Cela nous donne donc la liste de fonctions propres suivantes :
$$v_{\alpha_{1,1}}:= \left( x_1 +1 \right)^{1/2} \left(1-x_1-x_2 \right)^{1/2} \, e^{-\sqrt{-1} t_1},$$
$$v_{\alpha_{1,2}}:= \left( x_1 +1 \right)^{1/2} \left(x_2+1 \right)^{1/2} \, e^{-\sqrt{-1} (t_1 -t_2)},$$
$$v_{\alpha_{2,1}}:= \left( x_2 +1 \right)^{1/2} \left(1-x_1-x_2 \right)^{1/2} \, e^{-\sqrt{-1} t_2},$$
$$v_{\alpha_{2,2}}:= \left( x_2 +1 \right)^{1/2} \left(x_1+1 \right)^{1/2} \, e^{-\sqrt{-1} (t_2 -t_1)},$$
$$
v_{\alpha_{3,1}}:= \left( x_1 +1 \right)^{1/2} \left(1-x_1-x_2 \right)^{1/2} \, e^{\sqrt{-1} t_1},
$$
$$v_{\alpha_{3,2}}:= \left( x_2 +1 \right)^{1/2} \left(1-x_1-x_2 \right)^{1/2} \, e^{\sqrt{-1} t_2}.$$
\end{enumerate}
Remarquons que puisqu'il s'agit d'une métrique de Kähler-Einstein, ces fonctions propres correspondent aux fonctions propres du laplacien usuel comme calculées par exemple dans le chapitre 3 section C de \cite{berger1971spectre} mais exprimées dans le système de coordonnées actions-angles.

\subsection{$\CC \PP^2$ éclaté en un point}
On considère le cas où la variété complexe $M$ est  $\CC \PP^2$ éclaté en un point, que nous noterons par la suite $M_1$. On sait (voir chapitre 1 de \cite{guillemin1994moment})  que le polytope $P$ associé à cette variété torique est le trapèze $\tilde{\tau}$ donné par :
$$
\tilde{\tau} = \left\lbrace x \in \RR^2 ~:~ \forall i=1,\cdots,4~~ \langle x , \nu_i \rangle \geq -1 \right \rbrace,
$$
où
$$
\nu_1 := e_2 
,~~
\nu_2 := -e_1,~~
\nu_3 := e_1,~~
\nu_4 := e_1 -e_2.
$$
On peut déjà calculer les racines de Demazure de ce trapèze :
$$
\alpha_{1,1}:=e_2,~~ \alpha_{2,1}:=-e_1,~~ \alpha_{2,2}:=-e_1-e_2,~~ \alpha_{4,1}:=-e_2.
$$

De plus, on devra considérer le  trapèze $\tau$ qui se trouve être le trapèze translaté du trapèze $\tilde{\tau}$ par le vecteur $(2,1)$. On montre alors que $\sigma : (x,y) \in [1,3] \times [0,1] \mapsto (x,xy) \in \tilde{\tau}$ est un difféomorphisme. Cela entraîne que $M_1$ est une \textit{variété torique de Calabi}. Avant de définir cette notion, nous avons besoin de la définition d'un \textit{trapèze de Calabi}.

\bdefi

Un trapèze de Calabi est un polytope de $\RR^2$ qui est égal à l'image d'un rectangle $[\alpha_1, \alpha_2] \times [\beta_1, \beta_2] \subset \RR^2$ avec $\alpha_1>0$ et $\beta_1 \geq 0$ par l'application $\sigma : (x,y) \mapsto (x,xy)$.
\edefi

On fixe maintenant un trapèze de Calabi donné par l'image d'un rectangle $[\alpha_1, \alpha_2] \times [\beta_1, \beta_2]$. On peut alors écrire les normales à chaque arête sous la forme
$$
u_{\alpha_1} = C_{\alpha_1} \left( \begin{array}{c}
\alpha_1 \\ 
0
\end{array}  \right),~~ C_{\alpha_1}>0
$$
$$
u_{\alpha_2} = C_{\alpha_2} \left( \begin{array}{c}
\alpha_2 \\ 
0
\end{array}  \right),~~ C_{\alpha_2}<0
$$
$$
u_{\alpha_1} = C_{\beta_1} \left( \begin{array}{c}
\beta_1 \\ 
-1
\end{array}  \right),~~ C_{\beta_1}<0
$$
$$
u_{\beta_2} = C_{\beta_2} \left( \begin{array}{c}
\beta_2 \\ 
-1
\end{array}  \right),~~ C_{\beta_2}>0.
$$
Ainsi le couple formé du trapèze de Calabi et ses vecteurs normaux (appelé trapèze de Calabi étiqueté dans l'article \cite{Eve2}) est entièrement déterminé par 8 paramètres :
$$
(\alpha_1, \alpha_2, \beta_1, \beta_2, C_{\alpha_1}, C_{\alpha_2}, C_{\beta_1}, C_{\beta_2}).
$$

Si on prend le trapèze de Calabi dont les paramètres sont $(1,3,0,1,1,-\frac{1}{3}, -1,1)
$, alors on retrouve le polytope $\tilde{\tau}$ associé à $M_1$. On notera $\mu$ l'application moment telle que $\im \mu = \tau$ et $\tilde{\mu}$ l'application moment telle que $\im \tilde{\mu} = \tilde{\tau}$. En particulier, on a $\tilde{\mu} = \mu - (2,1)$.

Terminons en rappelant le résultat fondamental qui nous servira par la suite à chercher les solitons de Kähler-Ricci.

\bdefi
Soit $(M,\omega,J,\TT^2, \mu)$ une variété torique kählérienne compacte connexe de dimension réelle $4$. On dit que $M$ est torique de Calabi s'il existe deux fonctions $\TT^2$-invariantes $x$ et $y \in \C^{\infty}(M)$ strictement positives telles que l'application moment $\mu$ soit donnée par $\mu=(x,xy)$. On appelle alors $(x,y)$ les coordonnées de Calabi.
\edefi

On remarque que le polytope associé à $M$ est alors un trapèze de Calabi. Nous avons la proposition $4.9$ de \cite{Eve2}.
\bprop
Soit $(M,\omega,J,\TT^2, \mu)$ une variété torique de Calabi ayant pour polytope le trapèze de Calabi donné par les paramètres $(\alpha_1, \alpha_2, \beta_1, \beta_2, C_{\alpha_1}, C_{\alpha_2}, C_{\beta_1}, C_{\beta_2})$ et ayant pour coordonnées de Calabi $(x,y)$. Remarquons que $\Im \, x = [\alpha_1, \alpha_2]$ et $\Im \, y = [\beta_1, \beta_2]$. On pose $(t,s)$ les deux cordonnées d'angles correspondantes sur $M^{0}$. Alors nous avons qu'il existe deux fonctions $A \in \C^{\infty}(M)$ et $B \in \C^{\infty}(M)$ telles que
\begin{itemize}
\item[$\bullet$] les fonctions $A$ et $B$ ne dépendent respectivement que de la variable $x$ et $y$ et sont strictement positives sur $M^{0}$,
\item[$\bullet$] la métrique $g$, dite métrique de Calabi, a une expression particulière sur $M^{0}$ :
\begin{equation}\label{metpart}
g \vert_{M^{0}} = \cfrac{x}{A(x)} \, dx^2 +  \cfrac{x}{B(y)} \, dy^2 + \cfrac{A(x)}{x} \, (dt+y ds)^2 + x B(y) ds^2,
\end{equation}
\item[$\bullet$]  on a les conditions au bord suivantes pour $i=1,2$ :
\begin{equation}
A(\alpha_i)=0, ~~ B(\beta_i)=0,
\end{equation}
et
\begin{equation}
A'(\alpha_i)= \cfrac{2}{C_{\alpha_i}}, ~~ B'(\beta_i)= \cfrac{2}{C_{\beta_i}}.
\end{equation}
\end{itemize}
\cqfd
\eprop
De plus, on montre alors que
$$
Scal_g = \cfrac{A''(x)+B''(y)}{x}
$$
et
$$
\overline{Scal}:= \cfrac{4}{\alpha_1 + \alpha_2} \left( \cfrac{1}{\alpha_2 - \alpha_1} \left( \cfrac{1}{C_{\alpha_1}} - \cfrac{1}{C_{\alpha_2}} \right) - \cfrac{1}{\beta_2 - \beta_1} \left( \cfrac{1}{C_{\beta_1}} - \cfrac{1}{C_{\beta_2}} \right) \right).
$$
La recherche de solitons de Kähler-Ricci à l'aide du théorème précédent est faite dans l'article \cite{Eve1} où il est montré le résultat suivant (lemme $3.8$).
\blemm
Supposons que $g$ est une métrique torique  de Calabi pour $A \in \C^{\infty}([\alpha_1, \alpha_2])$ et $B \in \C^{\infty}([\beta_1, \beta_2])$. Pour tout $a \in \ttt$, en écrivant $\langle \mu,a \rangle= a_1 x + a_2 xy$, le couple $(g,a)$ est un soliton de Kähler-Ricci si
\begin{equation}\label{equasol}
-A''(x) -2 a_1 A'(x) - x \, \overline{Scal} = m , ~~ B''(y)=m,
\end{equation}
où $m=\cfrac{1}{\beta_2 - \beta_1} \left(  \cfrac{2}{C_{\beta_1}} - \cfrac{2}{ C_{\beta_2}} \right)$.
En particulier, on a que $a_2=0$.
\elemm
\rem Nos conventions sont différentes d'où la différence de signes dans l'équation \eqref{equasol}. \\

Revenons à $\CC \PP^2$ éclaté en 1 point. On sait qu'il n'admet pas de métriques de Kähler-Einstein (voir \cite{sze}) donc on doit chercher à résoudre l'équation \eqref{equasol} avec $a_1 \neq 0$. Rappelons que les paramètres de Calabi sont donnés par
$$
(1,3,0,1,1,-\cfrac{1}{3}, -1,1).
$$
Un calcul direct nous donne que
\begin{equation}\label{mpart}
m=-4,
\end{equation}
et
\begin{equation}\label{Spart}
\overline{Scal}=-4.
\end{equation}
On peut alors résoudre l'équation \eqref{equasol} et on obtient :
\begin{equation}
A(x)= \cfrac{-1}{{a_1}^3} \left[ \left( {a_1}^2 - \cfrac{1}{2} \right) \, \exp(-2{a_1} \, ( x -1) ) +{a_1}^2 \, x^2  - (2{a_1}^2 +{a_1}) \, x + \left( {a_1} +\cfrac{1}{2} \right) \right],
\end{equation}
et la fonction $B$ a pour expression
\begin{equation}
B(y)= -2 \, (y-1) \, y = -2 y^2 +2y.
\end{equation}
On peut aussi calculer $a_1$ grâce à l'invariant de Futaki et on obtient que $a_1$ est la solution non nulle de 
\begin{equation}
\left( a_1^2 - \cfrac{1}{2} \right) \, \exp(-4a_1) + 3 \, a_1^2 - 2 \,a_1 + \cfrac{1}{2}=0.
\end{equation}
On peut alors calculer la matrice $H$ associée au potentiel symplectique du soliton. Pour cela, on rappelle que nous avons (voir le lemme $3.8$ de \cite{Eve1})  
\begin{equation*}
H(x,y)
=  \cfrac{1}{x}
\left(
\begin{array}{cc}
A(x) & y \, A(x) \\ 
y \, A(x) & x^2 B(y) + y^2 \, A(x),
\end{array} 
\right)
\end{equation*}
Cette fonction $H$ peut aussi s'écrire sous la forme
\begin{equation}\label{Hequation}
H(\mu_1,\mu_2)
=  \cfrac{1}{\mu_1}
\left(
\begin{array}{ccc}
A(\mu_1) &~~~~ & \cfrac{\mu_2}{\mu_1} \, A(\mu_1) \\ 
\cfrac{\mu_2}{\mu_1} \, A(\mu_1)& ~~~~ & \mu_1^2 B(\mu_1,\mu_2) + \left( \cfrac{\mu_2}{\mu_1} \right)^2 \, A(\mu_1),
\end{array} 
\right)
\end{equation}
On calcule alors la matrice $G$ inverse de la matrice $H$ :
$$
G
=
\left(
\begin{array}{ccc}
\cfrac{x}{A(x)} + \cfrac{y^2}{xB(y)} & ~~ & -\cfrac{y}{x B(y)} \\ 
& &\\
-\cfrac{y}{x B(y)} & ~~ & \cfrac{1}{xB(y)} \\ 
\end{array} 
\right)
=
\left(
\begin{array}{ccc}
\cfrac{\mu_1}{A(\mu_1)} + \cfrac{\mu^2_2}{\mu_1^3B(\mu_1,\mu_2)} & ~~ & -\cfrac{\mu_2}{\mu_1^2 B(\mu_1,\mu_2)} \\ 
& & \\
-\cfrac{\mu_2}{\mu_1^2 B(\mu_1,\mu_2)} & ~~ & \cfrac{1}{ \mu_1 B(\mu_1,\mu_2)} \\ 
\end{array}
\right),
$$
où 
$$
A(\mu_1)= \cfrac{-1}{a^3} \left[ \left( a^2 - \cfrac{1}{2} \right) \, \exp(-2a \, ( \mu_1 -1) ) +a^2 \, \mu_1^2  - (2a^2 +a) \, \mu_1 + \left( a + \cfrac{1}{2} \right) \right]
$$
et
$$
B(\mu_1,\mu_2)= -2\left( \cfrac{\mu_2}{\mu_1} \right)^2 +2 \left( \cfrac{\mu_2}{\mu_1} \right).
$$
On peut alors se ramener au polytope algébrique en utilisant l'expression
$$
\mu = \tilde{\mu} + (2,1).
$$
Notons $\phi$ le potentiel symplectique associé à la métrique de Calabi telle que la matrice $G$ ci-dessus soit la matrice hessienne de $\phi$. Nous obtenons que les fonctions propres sont données par
$$v_{\alpha_{1,1}}:= \left( \tilde{\mu}_2 +1 \right) \cdot e^{(\nabla \phi)_2 + \sqrt{-1} \, t_2},$$
 $$v_{\alpha_{2,1}}:= \left(  - \tilde{\mu}_1 +1 \right) \cdot e^{(\nabla \phi)_1 - \sqrt{-1} \, t_1},$$
$$v_{\alpha_{2,2}}:= \left(  - \tilde{\mu}_1 - \tilde{\mu}_2 +1 \right) \cdot e^{(\nabla \phi)_1 + (\nabla \phi)_2 - \sqrt{-1} \, \left( t_1 + t_2 \right)},$$
$$v_{\alpha_{4,1}}:= \left( - \tilde{\mu}_2 +1 \right) \cdot e^{(\nabla \phi)_2 - \sqrt{-1} \, t_2}.$$

\bibliographystyle{alpha}
\bibliography{biblio}

\end{document}